\documentclass[preprint,authoryear]{elsarticle}

\biboptions{sort&compress}

\usepackage{graphicx}      % include this line if your document contains figures
\usepackage{natbib}        % required for bibliography
\usepackage{amssymb}
\usepackage{amsmath}
\usepackage{enumerate}
\usepackage{placeins}
\usepackage{hyphenat}
\usepackage{epstopdf}

\newtheorem{thm}{Theorem}
\newtheorem{lem}{Lemma}
\newtheorem{define}{Definition}
\newtheorem{corol}{Corollary}
\newproof{pf}{Proof}

\begin{document}
\bibliographystyle{model5-names}
\pagestyle{plain}
\setcounter{page}{0}
\pagenumbering{arabic}
\begin{frontmatter}

\title{Sufficient Conditions for Feasibility and Optimality of Real-Time Optimization Schemes -- I. Theoretical Foundations} 

\author[epfl]{Gene A. Bunin}
\ead{gene.bunin@epfl.ch}
\author[epfl]{Gr\'egory Fran\c cois}
\ead{gregory.francois@epfl.ch}
\author[epfl]{Dominique Bonvin\corref{cor1}}
\ead{dominique.bonvin@epfl.ch}

\cortext[cor1] {Author to whom correspondence should be addressed (tel.: +41 21 6933843, fax: +41 21 6932574).}

\tnotetext[abb]{Nonstandard abbreviations used in this article: SCFO (sufficient conditions for feasibility and optimality).}

\address[epfl]{Laboratoire d'Automatique, Ecole Polytechnique F\'ed\'erale de Lausanne, CH-1015 Lausanne, Switzerland}

\begin{abstract}                % Abstract of not more than 250 words.

The idea of iterative process optimization based on collected output measurements, or ``real-time optimization'' (RTO), has gained much prominence in recent decades, with many RTO algorithms being proposed, researched, and developed. While the essential goal of these schemes is to drive the process to its true optimal conditions without violating any safety-critical, or ``hard'', constraints, no generalized, unified approach for \emph{guaranteeing} this behavior exists. In this two-part paper, we propose an implementable set of conditions that can enforce these properties for any RTO algorithm. The first part of the work is dedicated to the theory behind the sufficient conditions for feasibility and optimality (SCFO), together with their basic implementation strategy. RTO algorithms enforcing the SCFO are shown to perform as desired in several numerical examples -- allowing for feasible-side convergence to the \emph{plant} optimum where algorithms not enforcing the conditions would fail.

\;


\end{abstract}

\end{frontmatter}
%===============================================================================

\section{Real-Time Optimization: Problem Structure, Characteristics, and Challenges}

The idea of optimization is present in many experimental settings, as it is often the case that a user/operator is interested in finding the ``best'' set of decision/design variables so as to minimize (maximize) a certain cost (profit) criterion. Noting that many such problems will also involve constraints that limit the decision/design space of the variables, we may state the resulting problem mathematically as follows:

\begin{equation}\label{eq:realopt}
\begin{array}{l}
\mathop {{\rm{minimize}}}\limits_{\bf{u}}\hspace{4mm}\phi_p ({\bf{u}}) \\
{\rm{subject}}\hspace{1mm}{\rm{to}}\hspace{3mm}{\bf{G}}_p({\bf{u}}) \preceq {\bf{0}} \\
\hspace{18.3mm}{\bf{u}}^L \preceq {\bf u} \preceq {\bf u}^U
\end{array}\hspace{3mm},
\end{equation}

\noindent with ${\bf u} \in \mathbb{R}^{n_u}$ denoting the manipulated decision variables (the ``inputs''), $\phi_p : \mathbb{R}^{n_u} \rightarrow \mathbb{R}$ the cost function to be minimized, ${\bf G}_p$ a set of $n_g$ inequality constraint functions $g_p : \mathbb{R}^{n_u} \rightarrow \mathbb{R}$, and ${\bf u}^L$ and ${\bf u}^U$ a set of lower and upper limits on the inputs (the box constraints). We adopt the subscript $p$ (i.e. the ``plant'') to signify that the corresponding functions are unknown, black-box, or uncertain, i.e. that one may only evaluate them at various discrete instants by applying various choices of ${\bf u}$, with each choice amounting to a single experiment.

Given the generality of (\ref{eq:realopt}), it is not surprising that a great number of applications give rise to this problem. Examples include:

\begin{itemize}
\item steady-state optimization \citep{Chen:87,Fatora:92,Naysmith:95,Cheng:04,Brdys:05,Gao:05,Flemming:07,Engell:07,Tatjewski:08}, where a plant operator is interested in finding operating conditions that maximize the profit of a dynamic plant at steady state,
\item optimization of a dynamic profile in a batch process \citep{Srinivasan:03,Srinivasan:03a,Francois:05,Kadam:07,Georgakis:09,Costello:11}, where the inputs are the ``handles'' used to parameterize, in some parsimonious fashion, the trajectory being optimized,
\item any application where experiments are carried out to optimize a response for which a model is not available, typically done using design-of-experiments and response-surface methods (see, e.g., \cite{Myers:09} and the examples therein),
\item controller tuning/design \citep{Hjarmarsson:98,Karimi:04,Krstic:06,Magni:09,Bunin:12a,Bunin:13}, where one is interested in finding the controller parameters that yield the best closed-loop performance,
\item numerical optimization where function evaluations are ``expensive'' \citep{Conn:09}, e.g. problems where evaluating a function is time-consuming as it involves simulating a system of differential equations \citep{Vugrin:03},
\end{itemize}

\noindent among others.

\begin{figure}
\begin{center}
\includegraphics[width=6cm]{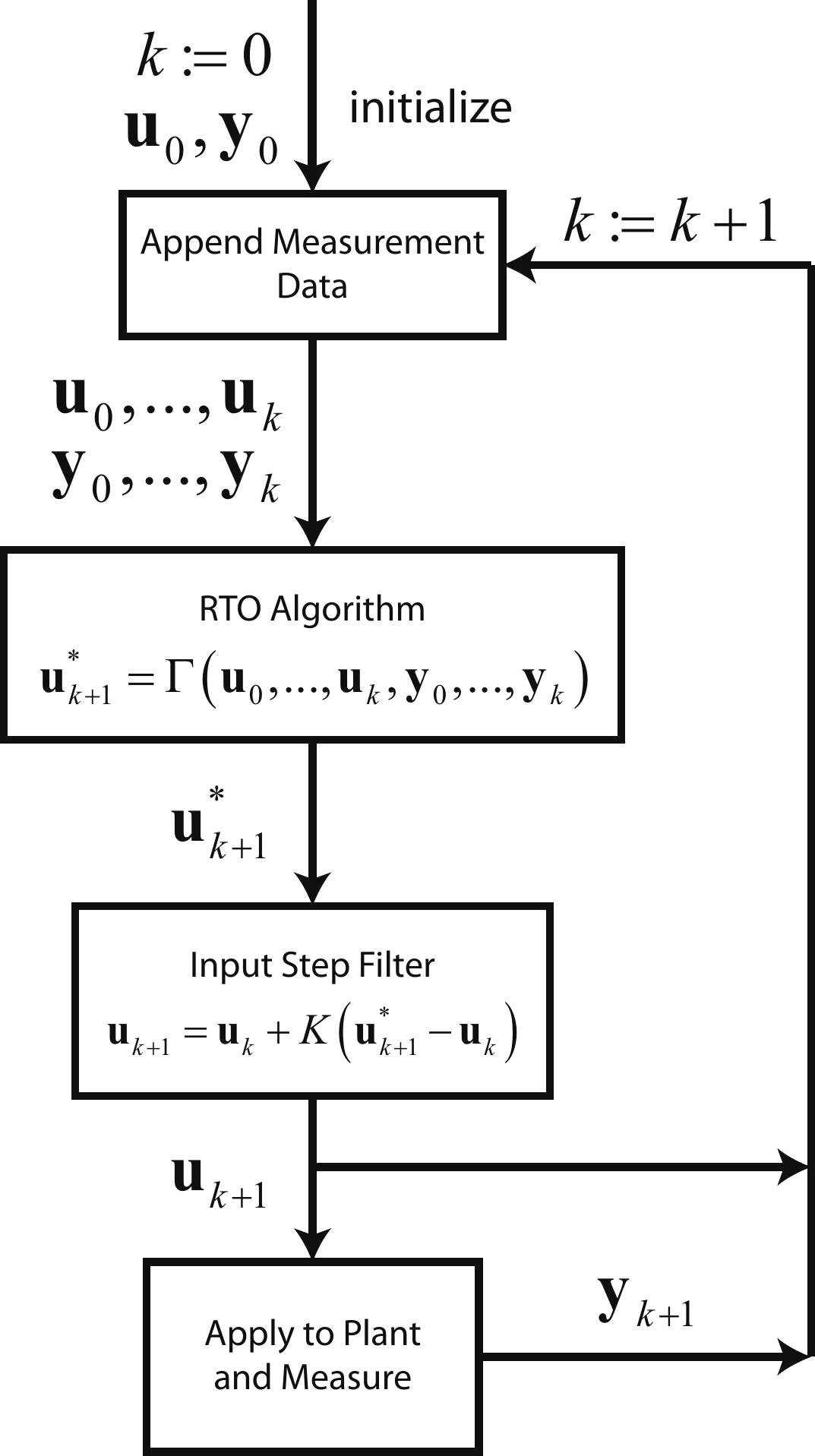}    % The printed column width  is 8.4 cm.
\caption{Generalized feedback structure of RTO schemes.}
\label{fig:fback}
\end{center}
\end{figure}

Because the functions involved are unknown, it is evident that one cannot hope to solve Problem (\ref{eq:realopt}) directly and without any experimentation. We may, however, attempt to solve the problem \emph{iteratively} -- that is, we may apply certain choices of ${\bf u}$, measure the results, and then apply ``more intelligent'' choices of ${\bf u}$ based on what we have learned. A general formulation of such an approach may be given as:

\begin{equation}\label{eq:algostep}
\begin{array}{l}
{\bf u}_{k+1}^* = \Gamma ({\bf u}_0,...,{\bf u}_k,{\bf y}_0,...,{\bf y}_k) \\
{\bf u}_{k+1} = {\bf u}_k + K ({\bf u}_{k+1}^* - {\bf u}_k)
\end{array}\hspace{3mm},
\end{equation}

\noindent where ${\bf y} \in \mathbb{R}^{n_y}$ denotes the measurements (the ``outputs''), $k$ the iteration (experiment) counter, and $\Gamma(\cdot)$ some prescribed adaptation law to yield an ``optimal'' target ${\bf u}_{k+1}^*$ for the next experiment, which is often filtered with a gain of $K \in [0,1]$ as a safety precaution \citep{Brdys:05}. It is precisely this iterative nature of (\ref{eq:algostep}) that introduces the ``real-time'' element. We will, hereafter, refer to $\Gamma(\cdot)$ as a \emph{real-time optimization (RTO) algorithm}, and to any problem having the form of (\ref{eq:realopt}) that is solved by (\ref{eq:algostep}) as an \emph{RTO problem}. A graphical illustration of the iterative RTO procedure is presented in Figure \ref{fig:fback}.

Unfortunately, RTO problems rarely have ``neat'' solutions, and a number of complications arise due to the experimental nature of the problem. The major culprit is, of course, the black-box nature of the problem -- one does not know the analytical properties of $\phi_p$ and ${\bf G}_p$, thereby making it impossible to optimize via even the simplest numerical approaches like the gradient descent, let alone via more elegant (e.g. interior-point) methods. A typical strategy in many RTO algorithms \citep{Tadej:01,Brdys:05,Gao:05,Marchetti:10,Rodger:10,Bunin:12} is to estimate the gradient from discrete measurements, but this is often made difficult by the presence of \emph{measurement noise}, thereby requiring a number of precautions and/or limitations for such an estimation to be effective. Another pressing difficulty is feasibility, since the inequality constraints ${\bf G}_p$ may be safety-critical or ``hard'', in which case one is strictly forbidden from applying any ${\bf u}$ that violates one of the inequalities, as doing so may promote hazardous conditions that either endanger the personnel or irreversibly damage the experimental equipment. Finally, the innately costly nature of performing experiments generally demands that RTO algorithms move quickly and that they move in the right direction -- an RTO algorithm that does not lead to immediate improvement may be abandoned after only a few iterations, while an algorithm that promises optimality, but only after 10,000 iterations, may never be implemented in the first place.

In spite of these challenges, a number of RTO algorithms have been proposed over the years and may be summarized (more or less completely) as follows:

\begin{itemize}
\item model-based approaches where a parametric model of Problem (\ref{eq:realopt}) is available and where the measurements are used for parameter estimation, after which the updated parametric model is optimized numerically to yield ${\bf u}_{k+1}^*$ \citep{Jang:87,Chen:87,Fatora:92,Naysmith:95},
\item model-based approaches where bias correction terms are estimated and added to the model, with a numerical optimization of the corrected model used to yield ${\bf u}_{k+1}^*$ \citep{Roberts:78,Brdys:05,Gao:05,Xu:08,Marchetti:09,Francois2013},
\item response-surface methods that are based on optimal experimental design \citep{Georgakis:09,Myers:09}, where a (usually quadratic) model is built from the collected measurements and then optimized (${\bf u}_{k+1}^*$ often becoming fixed once this is done),
\item direct-search methods that attempt to estimate and use the gradient \citep{Box:51,Garcia:84,Francois:05,Krstic:06,Bunin:12a}, akin to the standard gradient-descent algorithm in numerical optimization (${\bf u}_{k+1}^*$ being a step in the gradient descent direction),
\item derivative-free methods that avoid estimating the gradient entirely and optimize purely via zero-order knowledge \citep{Box:69,Alexandrov:97,Holmes:03,Conn:09}, with ${\bf u}_{k+1}^*$ defined differently depending on the algorithm used,
\end{itemize}

\noindent together with any hybrids/variations of the above.

However, even with this variety of methods, the state-of-the-art in RTO remains largely \emph{ad hoc}, as none of the above approaches are capable of overcoming the stated challenges and \emph{consistently} solving Problem (\ref{eq:realopt}) to local optimality without violating the hard constraints in the process. As such, while one algorithm may perform well for a certain problem, it may perform quite poorly (or even be dangerous) for another. The result is an additional burden on the user, as even \emph{conceptual} guarantees of feasibility and optimality are absent in these algorithms.

Our goal in this work is not to propose an algorithm that solves all of these issues -- rather, we aim to propose a set of sufficient conditions that serve to ensure the desired properties while being applicable to any RTO algorithm that falls into the framework of Figure \ref{fig:fback}. In doing so, we hope to help lay the foundations for an RTO \emph{theory} in a field that has often relied on heuristic approaches without rigorous guarantees. For convenience, we will refer to the full set of these conditions by the acronym SCFO (the ``sufficient conditions for feasibility and optimality'').

The basic idea of the SCFO is to ensure that the next applied iterate, ${\bf u}_{k+1}$, always fall into the local feasible descent space\footnote{The space where every choice of ${\bf u}$ is guaranteed to both decrease the cost and to satisfy the plant constraints locally.} and that the distance between ${\bf u}_{k}$ and ${\bf u}_{k+1}$ be small enough so as to ensure that the applied iterate will both be feasible and have an improved cost. The main contribution of this work is the theory and rigor behind this strategy (detailed in Theorems 2-5). The secondary contribution is the algorithm proposed for the basic SCFO implementation, which takes the input target provided by the RTO algorithm and projects it onto the local feasible descent space in such a manner so as to ultimately enforce the desired feasible-side convergence without introducing major performance drawbacks (i.e. very slow convergence speed).

In order to properly establish the aforementioned guarantees, we are forced to make several basic assumptions throughout this paper:

\begin{enumerate}[]
\item {\bf{A1}}: The functions $\phi_p$ and ${\bf{G}}_p$ are twice continuously differentiable over the \emph{relevant input space} $\mathcal{I} = \{ {\bf u} : {\bf u}^L \preceq {\bf u} \preceq {\bf u}^U \}$.
\item {\bf{A2}}: The initial point, ${\bf{u}}_0$, is strictly feasible with respect to the inequality constraints (${\bf{G}}_p ({\bf{u}}_0) \prec {\bf{0}}$).
\item {\bf{A3}}: For every function $g_{p,j}, j=1,...,n_g$, there exists an $\epsilon_{m,j} > 0$ such that, for any given iteration $k$, $0 \geq g_{p,j}({\bf{u}}_k) \geq -\epsilon_{m,j} \Rightarrow \nabla g_{p,j}({\bf{u}}_k) \neq {\bf{0}}$.
\item {\bf{A4}}: All RTO algorithms considered will never yield a target outside of $\mathcal{I}$, i.e. ${\bf u}_{k+1}^* \in \mathcal{I}$ always.
\end{enumerate}

\noindent Here, Assumption A1 is required for the existence of $1^{st}$- and $2^{nd}$-order upper bounds (see Section 2), as well as for the existence of the gradients of the cost and constraints, while Assumption A2 is required for feasibility (Section 3). Assumption A3 is required to formally avoid the (somewhat pathological) case where an active constraint has a gradient of zero, and is subsumed by the more standard linear independence constraint qualification \citep[e.g.][Ch. 9]{Fletcher:87}. Assumption A4 aids in simplifying some of the derivations and proofs in this work.

In this first part of the two-part contribution, the SCFO are stated and proven \emph{conceptually}. As such, more importance is given to theory and less to practical aspects. Particularly, the following are assumed to be available:

\begin{itemize}
\item Knowledge of the gradients of $\phi_p$ and ${\bf{G}}_p$ at the current iteration $k$.
\item Knowledge of global Lipschitz constants for ${\bf{G}}_p$.
\item Knowledge of a global quadratic upper bound for $\phi_p$.
\item Knowledge of the exact values of ${\bf{G}}_p$ at the current iteration $k$ (i.e. noise-free measurements, or perfect estimates, of the constraints).
\end{itemize}

\noindent The results presented in this first part are thus not directly applicable per se, and serve only to provide a theoretical base. Additionally, as our focus is solely on feasibility and optimality, the notion of convergence speed is also absent from our discussions. In the sequel paper \citep{Bunin:12c}, the practical implementation issues are treated in detail, considering the cases where the knowledge of the four items above is imperfect, as well as proposing modifications to the SCFO that may greatly speed up the convergence of a given algorithm.

The remainder of this paper is structured as follows. Following the presentation of the key mathematical concepts in Section 2, we dedicate Section 3 to the feasibility problem of meeting the hard constraints at every RTO iteration, which makes up the bulk of the paper and presents the core components of the theoretical contribution. Guaranteeing convergence to a Karush-Kuhn-Tucker (KKT) point via feasible iterates is then the subject of Section 4, and acts as a natural extension of the theory presented in the section before. Numerical examples are given in both sections to illustrate the methodology for a simple two-dimensional case with readily available geometric representations. Section 5 then serves to conclude the paper with a summary and discussion of the results and their potential significance in the RTO community.

\section{Mathematical Preliminaries}

Here, we present the tools and definitions that will act as the foundation for the results that follow in Sections 3 and 4.

A key necessity for the theory proposed in this work is to be able to upper bound the evolution (the worst-case possible growth) of twice continuously differentiable black-box functions, generalized here by $f : \mathbb{R}^{n_u} \rightarrow \mathbb{R}$, over the compact input space $\mathcal{I}$. For convenience and cohesion with results that come later, this is always done with respect to the input point ${\bf{u}}_k$, which may be thought of as the current iterate at the current RTO iteration $k$. The following two bounds, one linear and one quadratic, are proposed in Lemmas 1 and 2, respectively.

\begin{lem}{\bf (Linear Upper Bound on the Evolution of an Unknown Function)}

Let $f$ be continuously differentiable over $\mathcal{I}$, so that:

\begin{equation}\label{eq:lip0}
-\kappa_{i} < \frac{\partial f}{\partial u_i} \Big |_{\bf{u}} < \kappa_{i}, \hspace{3mm} \forall {\bf{u}} \in \mathcal{I}, \hspace{3mm} i = 1,...,n_u,
\end{equation}

\noindent where $\kappa$ are the univariate Lipschitz constants of $f$.

Then, the evolution of $f$ between any ${\bf{u}}_k, {\bf{u}}_{k+1} \in \mathcal{I} \setminus \left\{ {\bf{u}}_k, {\bf{u}}_{k+1} : {\bf{u}}_k = {\bf{u}}_{k+1} \right\}$ may be strictly bounded as:

\begin{equation}\label{eq:linbound}
f ({\bf{u}}_{k+1}) - f ({\bf{u}}_k) < \displaystyle \sum\limits_{i = 1}^{n_u} {\kappa_{i} | u_{k+1,i} - u_{k,i} | }.
\end{equation}

\end{lem}

\begin{pf}

Start by defining:

\begin{equation}\label{eq:Delta}
{\bf{\Delta}}_i = \left[ {\begin{array}{*{20}c}
   {\bf{I}}_i & {\bf{0}} \\
   {\bf{0}} & {\bf{0}}_{n_u-i} \\
\end{array}} \right],
\end{equation}

\noindent i.e. an $n_u \times n_u$ diagonal matrix with the first $i$ diagonal elements equal to 1 and the remaining $n_u - i$ equal to 0. This allows us to write the evolution of $f$ as the summation of its decoupled, univariate evolutions in the following compact form:

\begin{equation}\label{eq:unievol}
f ({\bf{u}}_{k+1}) - f ({\bf{u}}_k) = \displaystyle \sum\limits_{i = 1}^{n_u} \Big ( f({\bf{u}}_k + {\bf{\Delta}}_i ({\bf{u}}_{k+1} - {\bf{u}}_k)) - f({\bf{u}}_k + {\bf{\Delta}}_{i-1} ({\bf{u}}_{k+1} - {\bf{u}}_k)) \Big ). 
\end{equation}

From (\ref{eq:lip0}) and the definition of the Lipschitz constant for the univariate case, it follows that:

\begin{equation}\label{eq:unievolbound}
f({\bf{u}}_k + {\bf{\Delta}}_i ({\bf{u}}_{k+1} - {\bf{u}}_k)) - f({\bf{u}}_k + {\bf{\Delta}}_{i-1} ({\bf{u}}_{k+1} - {\bf{u}}_k)) \leq \kappa_{i} | u_{k+1,i} - u_{k,i} |,  
\end{equation}

\noindent with strict inequality when $u_{k+1,i} \neq u_{k,i}$. Since ${\bf u}_k \neq {\bf u}_{k+1}$, it follows that $\exists i : u_{k+1,i} \neq u_{k,i}$ and, summing (\ref{eq:unievolbound}) for $i= 1,...,n_u$, we immediately have the result in (\ref{eq:linbound}). \qed

\end{pf}

\begin{lem}{\bf (Quadratic Upper Bound on the Evolution of an Unknown Function)}

Let $f$ be twice continuously differentiable over $\mathcal{I}$, so that:

\begin{equation}\label{eq:M0}
-M_{ij} < \frac{\partial^2 f}{\partial u_i \partial u_j} \Big |_{\bf{u}} < M_{ij}, \hspace{3mm} \forall {\bf{u}} \in \mathcal{I}, \hspace{3mm} i,j = 1,...,n_u.
\end{equation}

Then, the evolution of $f$ between any ${\bf{u}}_k, {\bf{u}}_{k+1} \in \mathcal{I}$ may be bounded as:

\begin{equation}\label{eq:Qbound}
f({\bf{u}}_{k+1}) - f({\bf{u}}_{k}) \leq \nabla f({\bf{u}}_k)^T ({\bf{u}}_{k+1} - {\bf{u}}_{k}) + \frac{1}{2}({\bf{u}}_{k+1} - {\bf{u}}_{k})^T \overline {\bf{Q}} ({\bf{u}}_{k+1} - {\bf{u}}_{k}),
\end{equation}

\noindent where $\overline {\bf Q} \succ 0$ is the quadratic upper bound, a diagonal matrix with its diagonal elements defined as:

\begin{equation}\label{eq:Qbounddiag}
\overline Q_{ii} = \displaystyle \mathop {\sum} \limits_{j=1}^{n_u} M_{ij}, \;\; i = 1,...,n_u.
\end{equation}

\end{lem}

\begin{pf}

We refer the reader to the appendix. \qed

\end{pf}

Also common throughout the discussion that follows is the notion of descent halfspaces and approximately active ($\epsilon$-active) constraints. We define these as follows.

\begin{define}{\bf (Local Descent Halfspace)}

The strict local descent halfspace of $f$ at ${\bf{u}}_{k}$ is defined as the set $\{{\bf{u}} :$ $\nabla f({\bf{u}}_k)^T ({\bf{u}} - {\bf{u}}_{k}) < 0 \}$. Its nonstrict approximation is defined as $\{{\bf{u}} : \nabla f({\bf{u}}_k)^T$ $ ({\bf{u}} - {\bf{u}}_{k}) \leq -\delta, \; \delta > 0 \}$, with the quality of the approximation increasing as $\delta \rightarrow 0$.

\end{define}

\begin{define}{\bf($\epsilon$-Active Constraints)}

A constraint $g_{p,j}$ is said to be $\epsilon$-active at iteration $k$ if, for $\epsilon_j > 0$, $0 > g_{p,j}({\bf{u}}_k) \geq -\epsilon_j$. This is used to approximate an active constraint set by an $\epsilon$-active set, with the quality of approximation increasing as $\epsilon_j \rightarrow 0$.

\end{define}

We now combine the concepts of the descent halfspace and the quadratic upper bound to derive a sufficient step size and step direction to guarantee strict descent for $f$. This is crucial both to guarantee that an algorithm can preserve feasibility without converging prematurely to a constraint (Section 3) and to guarantee that an algorithm can decrease the cost at every iteration (Section 4).

\begin{lem}{\bf(Guaranteed Descent Step)}

Let ${\bf{u}}_{k+1}^*$ lie in the strict local descent halfspace of $f$ at ${\bf{u}}_{k}$ and let $K$ be the step size in the direction ${\bf{u}}_{k+1}^* - {\bf{u}}_{k}$, so that the next iterate is defined as in (\ref{eq:algostep}). A strict descent in the function value, i.e. $f({\bf{u}}_{k+1}) < f({\bf{u}}_{k})$, is guaranteed if the following condition holds:

\begin{equation}\label{eq:descentstep}
0 < K < -2\frac{\nabla f({\bf{u}}_k)^T ({\bf{u}}^*_{k+1} - {\bf{u}}_{k})}{({\bf{u}}^*_{k+1} - {\bf{u}}_{k})^T \overline {\bf{Q}} ({\bf{u}}^*_{k+1} - {\bf{u}}_{k})},
\end{equation}

\noindent with $\overline {\bf{Q}} \succ 0$ defined as in (\ref{eq:Qbounddiag}).

\end{lem}

\begin{pf}

The maximum value that $f({\bf{u}}_{k+1})$ can take is given by Lemma 2:

\begin{equation}\label{eq:Qboundw}
f({\bf{u}}_{k+1}) = f({\bf{u}}_{k}) + \nabla f({\bf{u}}_k)^T ({\bf{u}}_{k+1} - {\bf{u}}_{k}) + \frac{1}{2}({\bf{u}}_{k+1} - {\bf{u}}_{k})^T \overline {\bf{Q}} ({\bf{u}}_{k+1} - {\bf{u}}_{k}),
\end{equation}

\noindent where substituting the update law of (\ref{eq:algostep}) readily yields:

\begin{equation}\label{eq:Qboundw2}
f({\bf{u}}_{k+1}) = f({\bf{u}}_{k}) + K \nabla f({\bf{u}}_k)^T ({\bf{u}}_{k+1}^* - {\bf{u}}_{k}) + \frac{K^2}{2}({\bf{u}}_{k+1}^* - {\bf{u}}_{k})^T \overline {\bf{Q}} ({\bf{u}}_{k+1}^* - {\bf{u}}_{k}),
\end{equation}

\noindent which is quadratic in $K$ and satisfies $f({\bf{u}}_{k+1}) = f({\bf{u}}_{k})$ for:

\begin{equation}\label{eq:zeros}
K = 0, -2\frac{\nabla f({\bf{u}}_k)^T ({\bf{u}}^*_{k+1} - {\bf{u}}_{k})}{({\bf{u}}^*_{k+1} - {\bf{u}}_{k})^T \overline {\bf{Q}} ({\bf{u}}^*_{k+1} - {\bf{u}}_{k})},
\end{equation}

\noindent with the latter forced to be strictly positive from the definition of a strict local descent halfspace (i.e. strictly negative numerator) and the strict positivity of the denominator.

As (\ref{eq:Qboundw2}) must be strictly convex in $K$, it follows that $f({\bf{u}}_{k+1}) < f({\bf{u}}_{k})$ for all values of $K$ between the two bounds of (\ref{eq:zeros}). \qed

\end{pf}

Finally, the algorithm that enforces the SCFO relies on projecting, at every iteration, the RTO target onto a local feasible descent space, which is characterized by an intersection of halfspaces. To prove that this projection will always be feasible unless a KKT point has been reached, the following theorem of alternatives will be needed.

\begin{thm}{\bf(Gale's Theorem)}

Let ${\bf{J}}_k = \left[ \nabla f_1({\bf{u}}_k) ... \nabla f_n({\bf{u}}_k) \right] $ and consider the set $\{{\bf{u}} : {\bf{J}}_k^T ({\bf{u}} - {\bf{u}}_{k}) \preceq -\boldsymbol{\delta} \}$. This set is empty and does not admit a solution iff:

\begin{equation}\label{eq:theoremalt}
\exists \boldsymbol{\nu} \in \mathbb{R}^{n}_+ : \hspace{1mm}\displaystyle \sum\limits_{i = 1}^{n} {\nu_{i}\nabla f_i({\bf{u}}_k)  } = {\bf{0}}, \hspace{3mm} -\boldsymbol{\nu}^T  \boldsymbol{\delta} < 0.
\end{equation}

\end{thm}

\begin{pf}

This is a well-known result in linear programming \citep[Th. 22.1]{Rockafellar:70}. \qed

\end{pf}

\section{Feasibility Guarantees}

The goal of this section is to focus on the guaranteed satisfaction of hard constraints, i.e. on constraints that, for safety reasons, must \emph{never} be violated\footnote{So as to work with the most limiting case, we will hereafter assume that all of the constraints in ${\bf G}_p$ are hard, i.e. it is required that ${\bf G}_p ({\bf u}_k) \preceq {\bf 0},\; \forall k$.}. We start with a brief review of the existing methods for enforcing RTO feasibility, and then proceed to propose a new feasibility-guaranteeing condition based on the use of an \emph{adaptive} filter gain in the input adaptation step (\ref{eq:algostep}), with the filter gain at iteration $k$ a function of the constraint values at iteration $k$ and the Lipschitz constants for those constraints. A drawback of the proposed approach is in its potentially premature convergence, and so we propose conditions that allow the algorithm to avoid this problem. Several examples are then given to illustrate the presented ideas. 

\subsection{Approaches to Guaranteeing Feasibility}

Some of the more theoretically elegant approaches to enforcing the feasibility of (\ref{eq:realopt}) include the stochastic (``chance constraint'') and worst-case robust approaches, which require a model of the constraints and are generally used in an initial design phase \citep{Beyer:07}, though they may be used for RTO needs when the model is adapted \citep{Zhang:02,Flemming:07}. In all cases, parametric uncertainty in the model is assumed and quantified. As a simple illustration of what this entails, suppose that a single constraint $g_{p,j}$ is uncertain due to a single uncertain parameter $\theta$, known to lie in the set $\Theta$, and that we would like to guarantee the following:

\begin{equation}\label{eq:congar0}
g_{p,j}({\bf{u}}_{k+1}) \leq 0,
\end{equation}

\noindent by guaranteeing:

\begin{equation}\label{eq:congar2}
\mathop{\sup}\limits_{\theta \in \Theta} g_j ({\bf{u}}_{k+1},\theta) \leq 0,
\end{equation}

\noindent where $g_j : \mathbb{R}^{n_u} \rightarrow \mathbb{R}$ is the parametric model.

The standard stochastic approach supposes a probability distribution for the possible values of $\theta$ and attempts to satisfy (\ref{eq:congar2}) with a specified probability, while the worst-case approach essentially attempts to do the same but with a probability of 1. The latter is generally considered as being very conservative since it accounts for all possible deviations, while the former suffers from the lack of robustness introduced by allowing the constraint to be violated, if only with a low probability.

While these methods are theoretically sound, there are several difficulties in applying them to RTO problems \citep{Quelhas:12}. The immediate issue that should come to mind is the initial assumption that the uncertainty is parametric, which will not hold in the general case. When structural uncertainties are present, it follows that both approaches lose their theoretical rigor, with the degree of the loss varying, naturally, with the amount of structural uncertainty. Additionally, even when the parametric uncertainty assumption is satisfied, obtaining the probability distributions or worst-case upper and lower bounds on the parameters may be a very challenging task. Finally, the problem may become computationally intractable, depending on the probability distribution and the way in which the parameters enter into the constraint model. This latter becomes particularly troublesome as the number of uncertain parameters and constraints grows. While it is sometimes proposed to solve this problem by linearizing the constraint model with respect to the uncertain parameters \citep{Zhang:02}, doing so may lead to significant inaccuracies when the process is nonlinear -- even if there is no structural uncertainty.

Some simpler methods of enforcing feasibility involve the use of back-offs on the constraints \citep{Govatsmack:05,Quelhas:12}:

\begin{equation}\label{eq:backoff}
g_j ({\bf{u}}_{k+1},\theta) + c \leq 0,\hspace{2mm} c > 0,
\end{equation}

\noindent or the use of limited input changes for each iteration (box constraints on the input step sizes \citep{Cheng:00,Gao:05}):

\begin{equation}\label{eq:boxstep}
{\bf{u}}_k - \Delta {\bf{u}}_{max} \preceq {\bf{u}}_{k+1} \preceq {\bf{u}}_k + \Delta {\bf{u}}_{max}.
\end{equation}

\noindent Both may be conservative, however, and are ultimately \emph{ad hoc} approaches that do not come with the desired guarantees of feasibility. 

In view of all this, \cite{Bunin:11b} proposed a simple but robust approach that uses a variable filter gain $K_k$ in the standard input filtering law (\ref{eq:algostep}):

\begin{equation}\label{eq:inputfilter}
{\bf{u}}_{k+1}  = {\bf{u}}_k  + K_k \left( {\bf{u}}^*_{k + 1} - {\bf{u}}_k \right),
\end{equation}

\noindent where the gain $K_k$ is defined at iteration $k$ in such a way so as to guarantee $g_{p,j}({\bf{u}}_{k+1}) < 0$ provided that $g_{p,j}({\bf{u}}_k) < 0$. The following theorem, largely taken from our previous work \citep{Bunin:11b}, derives the adaptive upper bound on the filter gain $K_k$ so as to guarantee recursive feasibility for any RTO scheme.

\begin{thm}{\bf (Sufficient Condition for Feasibility)}

Suppose that ${\bf{G}}_p ({\bf{u}}_k) \prec {\bf{0}}$, and that an RTO algorithm provides a new optimal input target ${\bf{u}}^*_{k + 1} \neq {\bf u}_k$ that is then filtered according to (\ref{eq:inputfilter}) to give ${\bf{u}}_{k + 1}$. Then, enforcing the following upper bound on the adaptive filter gain:

\begin{equation}\label{eq:gainupper}
K_{k} \le \mathop{\min}\limits_{j = 1,...,n_g} \left[ \frac{{-g_{p,j}({\bf{u}}_{k } )}}{\displaystyle \sum\limits_{i = 1}^{n_u} {\kappa_{ji} | u^*_{k+1,i} - u_{k,i} | }} \right],
\end{equation}

\noindent with $-\kappa_{ji} < \frac{\partial g_{p,j}}{\partial u_i} \Big |_{\bf u} < \kappa_{ji}, \forall {\bf u} \in \mathcal{I}$, is sufficient to guarantee that ${\bf{G}}_p ({\bf{u}}_{k+1}) \prec {\bf{0}}$.

\end{thm}

\begin{pf}

Noting that $K_k \in [0,1]$ by definition, we treat two cases: $K_k = 0$ and $K_k > 0$. For the former, we have the trivial result that ${\bf u}_{k+1} = {\bf u}_k \Rightarrow {\bf G}_p ({\bf u}_{k+1}) = {\bf G}_p ({\bf u}_{k}) \preceq {\bf 0}$, with feasibility guaranteed since the inputs are not adapted.

When $K_k > 0$ it follows that ${\bf u}_{k+1} \neq {\bf u}_k$, which allows us to use Lemma 1 to write:

\begin{equation}\label{eq:proofA2}
g_{p,j} ({\bf{u}}_{k+1}) < g_{p,j} ({\bf{u}}_{k}) + \displaystyle \sum\limits_{i = 1}^{n_u} {\kappa_{ji} | u_{k+1,i} - u_{k,i} | },\hspace{3mm} \forall j = 1,...,n_g.
\end{equation}

Substituting in the filter law (\ref{eq:inputfilter}) for ${\bf{u}}_{k+1}$ then leads to:

\begin{equation}\label{eq:proofA3}
g_{p,j} ({\bf{u}}_{k+1}) < g_{p,j} ({\bf{u}}_{k}) + K_k \displaystyle \sum\limits_{i = 1}^{n_u} {\kappa_{ji} | u^*_{k+1,i} - u_{k,i} | }, \hspace{3mm}\forall j = 1,...,n_g.
\end{equation}

Let us now choose $K_k$ in such a way so as to make the right-hand side less than or equal to 0, which will clearly guarantee the same (with strict inequality) for $g_{p,j} ({\bf{u}}_{k+1})$:

\begin{equation}\label{eq:proofA4}
g_{p,j} ({\bf{u}}_{k}) + K_k \displaystyle \sum\limits_{i = 1}^{n_u} {\kappa_{ji} | u^*_{k+1,i} - u_{k,i} | } \leq 0, \hspace{3mm}\forall j = 1,...,n_g.
\end{equation}

Rearranging leads to:

\begin{equation}\label{eq:proofA5}
K_k \leq \frac{-g_{p,j}({\bf{u}}_{k } )}{\displaystyle \sum\limits_{i = 1}^{n_u} {\kappa_{ji} | u^*_{k+1,i} - u_{k,i} | }}, \hspace{3mm}\forall j = 1,...,n_g.
\end{equation}

Now, to guarantee that this holds for all of the constraints, we simply take the component-wise minimum, which leads to the expression in (\ref{eq:gainupper}).\qed

\end{pf}

\vspace{2mm}
\noindent {\it{Remarks}}
\vspace{2mm}

\begin{figure}
\begin{center}
\includegraphics[width=7cm]{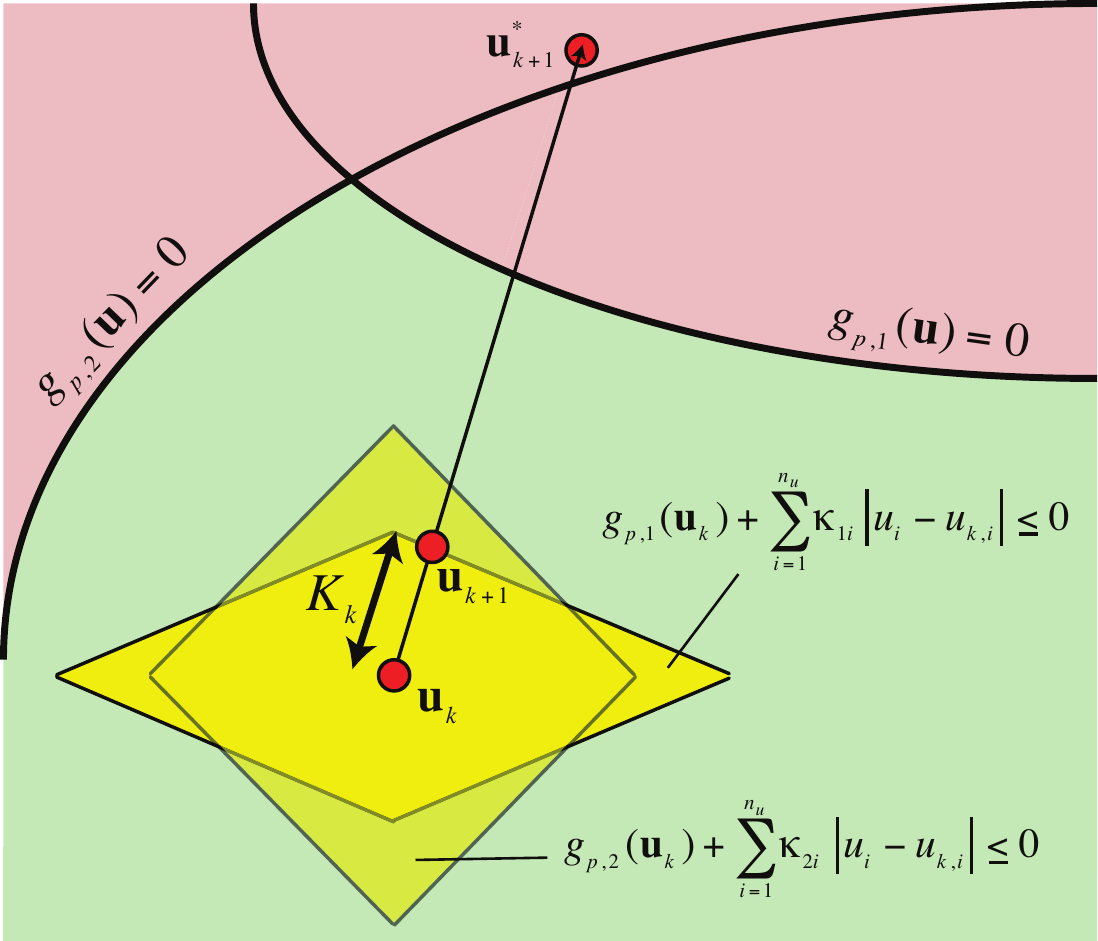}    % The printed column width  is 8.4 cm.
\caption{Geometric illustration of adaptive input filtering as a means of guaranteeing feasibility from iteration to iteration. Here, two constraints are considered. The green region denotes the feasible space, with the yellow polytopes corresponding to the robust feasible areas that are generated by the Lipschitz constants for each constraint. The adaptive filter gain value may be any value that keeps all of the constraint values in their worst-case growth regions. Here, the first constraint proves to be limiting.}
\label{fig:feasgraph1}
\end{center}
\end{figure}

\begin{itemize}
\item We can now justify the initial supposition of ${\bf{G}}_p ({\bf{u}}_k) \prec {\bf{0}}$ \emph{a posteriori}, since Assumption A2 calls for ${\bf{G}}_p ({\bf{u}}_0) \prec {\bf{0}}$ and recursive feasibility is enforced by (\ref{eq:gainupper}) up to iteration $k$ and beyond.
\item The main concept of Theorem 2 is very simple. Namely, we bound the worst-case evolution of the constraints (via the Lipschitz constants), calculate how far we can move before one of them, in the worst case, reaches its boundary, and then define the filter gain as the adaptive parameter with which to control this distance. The benefits of the approach lie in its generality, as it does not require the uncertainty to have any particular structure -- everything is simply lumped into the worst-case growth of the constraints. The adaptive element is also an advantage, in that the gain is always a function of the current measured/estimated values of the constraints. As such, when the plant is far away from the constraints, the gain is naturally larger, and is reduced as the constraints are approached. In the case of model-based RTO algorithms, this method also avoids introducing additional constraints into the original model-based optimization problem, thus guaranteeing that tractability is never lost. Figure \ref{fig:feasgraph1} illustrates the main idea of Theorem 2 geometrically.
\end{itemize}

We will, for the remainder of this section, adopt the following heuristic choice of $K_k$:

\begin{equation}\label{eq:gainuppereq}
K_{k} := \mathop{\min}\limits_{j = 1,...,n_g} \left[ \frac{{-g_{p,j}({\bf{u}}_{k } )}}{\displaystyle \sum\limits_{i = 1}^{n_u} {\kappa_{ji} | u^*_{k+1,i} - u_{k,i} | }} \right],
\end{equation}

\noindent as taking the largest allowable $K_k$ generally leads to faster progress for the RTO algorithm.

\subsection{Premature Convergence of the Adaptive Input Filter Method}

By themselves, the input filtering law (\ref{eq:inputfilter}) and the upper bound on the filter gain (\ref{eq:gainupper}) are sufficient for feasibility, provided that the initial point ${\bf{u}}_0$ is feasible. There is, however, a major algorithmic issue with using the inequality (\ref{eq:gainupper}), which may be illustrated as follows. Suppose that a single constraint $g_{p,j}$ is active and its value is arbitrarily close\footnote{We insist on ``arbitrarily close'' since the Lipschitz bounds are strict and so the constraint never reaches 0 exactly with the input filtering scheme as proposed in this work.} to 0. In looking at the expression for the adaptive filter gain in (\ref{eq:gainupper}), it becomes clear that the bound on $K_k$ will be arbitrarily close to 0 as well. This, in turn, implies that the input is not adapted and no further progress is made.

\begin{figure}
\begin{center}
\includegraphics[width=7cm]{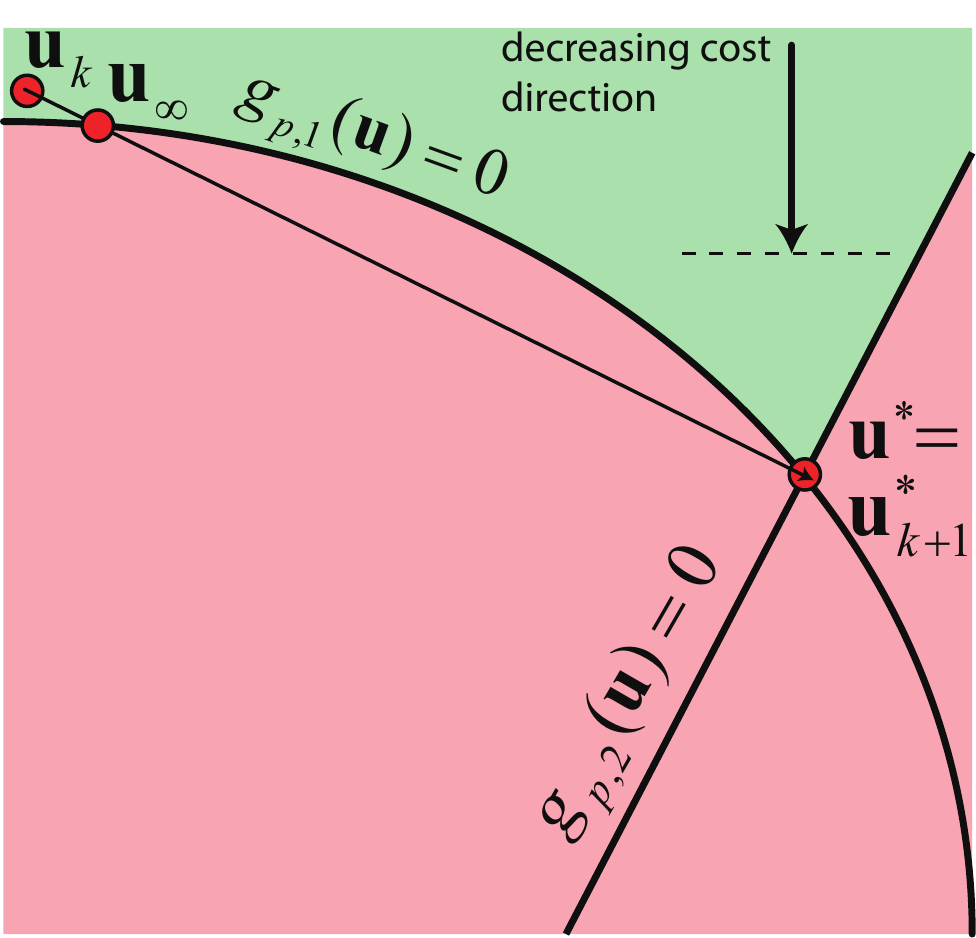}    % The printed column width  is 8.4 cm.
\caption{An example illustrating premature convergence, where one of the constraints becomes active and brings the filter gain to an arbitrarily small value. As the segment [${\bf u}_{\infty}$,${\bf u}_{k+1}^*$] is infeasible, the algorithm converges to ${\bf u}_{\infty}$ without being able to proceed further.}
\label{fig:feasbreak}
\end{center}
\end{figure}

This issue has been pointed out previously \citep{Bunin:11b}, and is illustrated in Figure \ref{fig:feasbreak} via a simple constructed example. Here, the RTO algorithm is perfect in that it provides the true plant optimum, ${\bf u}^*$, at each iteration. However, because of the concave nature of one of the plant constraints, it is clear that the iterates will eventually touch the boundary and remain arbitrarily close to it without being able to advance any further.

Depending on the particular problem, the losses incurred from such premature convergence may or may not be significant. In \cite{Bunin:11b}, it was argued that the majority of the optimality gains could be attained before one of the constraints became active, and this was shown through a numerical example. However, this may not hold always, and for processes that operate for a great number of iterations converging prematurely without any further improvement may not be an attractive result, even if feasibility is retained throughout. In the following section, we show how using the constraint \emph{gradients} allows us to construct an algorithm that avoids premature convergence by ``sliding'' along any constraints that are close to active.

\subsection{Avoiding Premature Convergence}

Prior to presenting an improved version of the method, it is first necessary to characterize the degree to which the Lipschitz constants of the constraint functions are strict.

\begin{define}{\bf (Degree of Strictness of the Lipschitz Constants)}

Denoting by $\tilde \kappa$ the nonstrict Lipschitz constants, which allows the nonstrict analogue to (\ref{eq:lip0}):

\begin{equation}\label{eq:lipNS}
- \tilde \kappa_{ji} \leq \frac{\partial g_{p,j}}{\partial u_i} \Big |_{\bf{u}} \leq \tilde \kappa_{ji}, \hspace{3mm} \forall {\bf{u}} \in \mathcal{I}, \hspace{3mm} i = 1,...,n_u, \;\; j = 1,...,n_g,
\end{equation}

\noindent the degree of strictness of the Lipschitz constants, $\gamma_\kappa$, is defined as:

\begin{equation}\label{eq:gammakappa0}
\gamma_\kappa = \mathop {\max} \limits_{\begin{array}{l}i=1,...,n_u\\ j = 1,...,n_g \end{array}} \frac{\tilde \kappa_{ji}}{\kappa_{ji}}.
\end{equation}

\end{define}

We now give the sufficient conditions for maintaining $K_k$ above a certain strictly positive value, thereby precluding the possibility of its becoming arbitrarily small, and thus avoiding premature convergence.

\begin{thm}{\bf (Sufficient Conditions for a Strictly Positive Adaptive Filter Gain)}

Let $K_k$ be chosen as in (\ref{eq:gainuppereq}). If ${\bf{u}}^*_{k+1}$ lies in the strict local descent halfspace for all of the $\epsilon$-active constraints:

\begin{equation}\label{eq:suff1}
\nabla g_{p,j}({\bf{u}}_k)^T ({\bf{u}}^*_{k+1} - {\bf{u}}_{k}) \leq -\delta_{g,j}, \hspace{1mm} \forall j : g_{p,j}({\bf{u}}_k) \geq -\epsilon_j,
\end{equation}

\noindent with $\delta_{g,j} > 0$, then:

\begin{equation}\label{eq:Klimit}
K_k > \mathop {\min} \limits_{j = 1,...,n_g} \left[ \frac{ \mathop {\min} \left( (1-\gamma_\kappa) \epsilon_j, (1-\gamma_\kappa) \displaystyle \frac{K_{\epsilon,j} \delta_{g,j}}{\gamma_\kappa}, -g_{p,j}({\bf{u}}_0 ) \right)  }{{\boldsymbol \kappa}_j^T ({\bf u}^U - {\bf u}^L)} \right] > 0,
\end{equation}

\noindent where:

\begin{equation}\label{eq:Keps}
K_{\epsilon,j} = \frac{2\delta_{g,j}}{({\bf{u}}^{U} - {\bf{u}}^{L})^T \overline {\bf{Q}}_j ({\bf{u}}^{U} - {\bf{u}}^{L}) },
\end{equation}

\begin{equation}\label{eq:kappavec}
{\boldsymbol \kappa}_j^T = \left[ \kappa_{j1} \; ... \; \kappa_{jn_u} \right],
\end{equation}

\noindent and $\overline {\bf{Q}}_j \succ 0$ denotes the quadratic upper bound for constraint $g_{p,j}$ as defined in Lemma 2. 

\end{thm}

\begin{pf}

The theorem may be proven by exploiting the result of Lemma 3 and showing that enforcing Condition (\ref{eq:suff1}) makes it impossible for any constraint, and thereby $K_k$, to approach 0 since for a small enough step ($K_{\epsilon,j}$) one will be forced to push the constraint away from activity. We refer the reader to the appendix for the detailed proof. \qed

\end{pf}

Several points of the theorem merit some discussion.

\vspace{2mm}
\noindent {\it{Geometric Interpretation of Theorem 3}}
\vspace{2mm}

\begin{figure}
\begin{center}
\includegraphics[width=5cm]{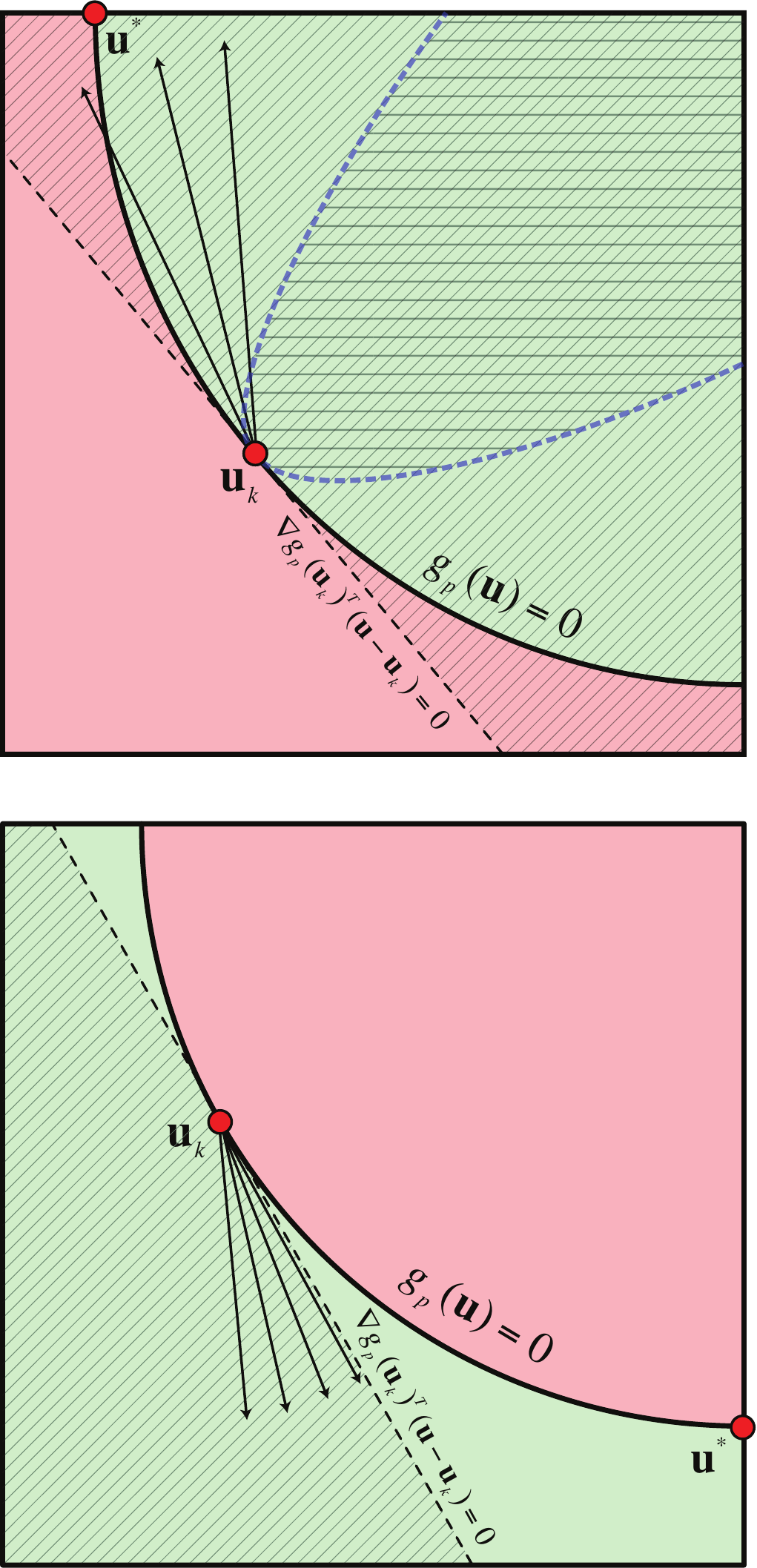}    % The printed column width  is 8.4 cm.
\caption{A geometric interpretation of the sufficient condition for avoiding premature convergence in the adaptive input filter scheme. In both examples, the lined area represents the local descent halfspace for the relevant constraint, with the arrows used to show potential directions that satisfy Condition (\ref{eq:suff1}) and for which sufficiently small step sizes will lead to true descent for the plant constraint. In the top case, a second lined area, corresponding to a hypothetical constraint with stronger nonlinear behavior, is also given, and it is clear how the same directions for this alternate case would require much smaller steps to guarantee true descent. The plant optimum, ${\bf u}^*$, is only given as a visual aid, and does not play any theoretical role here (as the focus is only on feasibility).}
\label{fig:feasgraph2}
\end{center}
\end{figure}

The geometric interpretation of the above theorem is presented via two examples in Figure {\ref{fig:feasgraph2}}, and illustrates the role of Condition (\ref{eq:suff1}) in forcing the optimal input direction to lie in the strict local descent halfspace of an $\epsilon$-active constraint. The top example in Figure {\ref{fig:feasgraph2}} shows this for a single convex constraint. When Condition (\ref{eq:suff1}) is met and the optimal direction lies in the strict local descent halfspace (the lined area), it is easy to see that each optimal direction will result, even if very locally, in a descent direction. The degree of locality here depends on both the direction (how orthogonal it is to the hyperplane $\nabla g_p({\bf{u}}_k)^T ({\bf{u}} - {\bf{u}}_{k}) = 0$) and the magnitude of the higher-order derivatives (in some sense, $\overline {\bf{Q}}$). Concerning the effect of direction, it is easy to see that directions orthogonal to the hyperplane would be steepest descent directions, with the constraint value going down substantially locally. With regard to the size of the nonlinear terms, we can imagine what would happen if, for example, the higher-order derivatives of the constraint were much larger -- this would squeeze the feasible space in the top example in Figure {\ref{fig:feasgraph2}} (shown via the dashed blue lines and the double-lined region) and lead to both much shorter feasible and descent steps if we were to follow the same optimal directions (arrowed lines). Both of these cases are reflected analytically in the upper bound (\ref{eq:descentstep}) of Lemma 3 -- steepest descent directions leading to a greater numerator and thus larger steps, and larger higher-order derivatives (greater $\overline {\bf{Q}}$) leading to a greater denominator and thus smaller steps. Finally, the illustration makes it obvious that any optimal direction that consistently points into the local \emph{ascent} halfspace would lead to the constraint value approaching 0, thereby forcing the algorithm to converge prematurely.

\vspace{2mm}
\noindent {\it{Significance of $\overline {\bf{Q}}$}}
\vspace{2mm}

It has been supposed throughout, both in the proofs of Lemma 3 and Theorem 3, that $\overline {\bf{Q}} \succ 0$, as this represents the worst case with respect to the evolution of $K_k$. If the constraint is strictly convex, then its higher-order evolution will be positive definite and such worst-case behavior will be justified (indeed, this is the case in the top example of Figure {\ref{fig:feasgraph2}}). We now consider the alternative when $\overline {\bf{Q}} \preceq 0$, which could be a valid upper bound when the constraint is concave (e.g. a negative semidefinite quadratic). If we allow such a bound, then it is clear from the analysis in the proof of Lemma 3 that $g_{p,j}({\bf{u}}_{k+1}) < g_{p,j}({\bf{u}}_k)$ will hold for \emph{all} strictly positive values $K_k$, and indeed it becomes clear why if we consider the interpretation in the bottom example of Figure {\ref{fig:feasgraph2}}. With a concave constraint, the hyperplane $\nabla g_p({\bf{u}}_k)^T ({\bf{u}} - {\bf{u}}_{k}) = 0$ supports the entire infeasible side of the constraint, and the condition in (\ref{eq:suff1}) thereby guarantees that any direction taken will only lead to both feasible and smaller values, regardless of the filter gain. The strict descent halfspace is global in this case.

More generally speaking, the abstract purpose of $\overline {\bf{Q}}$ is to help provide the conceptual proof that Condition (\ref{eq:suff1}) will prevent $K_k$ from going to 0 by always allowing a local descent direction whose length will depend, in large part, on $\overline {\bf{Q}}$. The fact that $\overline {\bf{Q}}$ may never be known in practice therefore does not hinder us -- we only need it to exist. As such, one does not need to go through the trouble of estimating $\overline {\bf{Q}}$ and attempting to enforce the step size derived in Lemma 3, although being able to do so would provide a tool to guarantee immediate descent at iteration $k$.

\vspace{2mm}
\noindent {\it{The Zero-Gradient Case}}
\vspace{2mm}

It is clear that the Condition (\ref{eq:suff1}) cannot possibly be satisfied if $\nabla g_{p,j}({\bf{u}}_k) $ $= {\bf{0}}$. Assumption A3 is used to avoid this case formally for a small enough choice of ${\boldsymbol \epsilon}$. 

\vspace{2mm}
\noindent {\it{The Strictness of the Lipschitz Constants}}
\vspace{2mm}

Not surprisingly, the strictness of the Lipschitz constants $\kappa$ is a \emph{necessity}, as the absence of a strict linear bound (as in Lemma 1) allows the possibility of a constraint becoming active, with its value exactly at 0. This is reflected in the bound on $K_k$ in (\ref{eq:Klimit}), which approaches 0 as $\gamma_\kappa \rightarrow 1$.

\vspace{2mm}

We next discuss how these ideas may be enforced algorithmically.

\subsection{Basic Implementation}

In general, there is no RTO algorithm that satisfies Condition (\ref{eq:suff1}) implicitly. However, the target optimum ${\bf{u}}^*_{k+1}$ may be seen as a degree of freedom that could be manipulated to fulfill (\ref{eq:suff1}) -- it is, in fact, the only degree of freedom, as the local gradient and current input are both fixed and cannot be modified. What we can do then is to project the optimum given by the RTO algorithm at iteration $k$ in the appropriate manner\footnote{The norm in the projection may be chosen depending on user preference and any necessary scaling that may need to be taken into account (i.e. the inputs should be adjusted so that they are of comparable magnitudes and have comparable sensitivities in the objective norm function). For all of the examples that follow in this article, we use a squared 2-norm as the default for any projections, with the chosen inputs already well scaled.}:

\begin{equation}\label{eq:feasproj}
\begin{array}{l}
\bar {\bf{u}}^*_{k+1} = {\rm{arg}} \mathop {{\rm{minimize}}}\limits_{{\bf{u}}}\hspace{4mm}\left\| {\bf{u}} - {\bf{u}}^*_{k+1} \right\|_2^2 \vspace{1mm}  \\
\hspace{17mm}{\rm{subject}}\hspace{1mm}{\rm{to}}\hspace{3mm}\nabla g_{p,j}({\bf{u}}_k)^T ({\bf{u}} - {\bf{u}}_{k}) \leq -\delta_{g,j} \vspace{1mm} \\
\hspace{36mm}\forall j: g_{p,j}({\bf{u}}_k) \geq -\epsilon_j \vspace{1mm}\\
\hspace{36mm}{\bf{u}}^L \preceq {\bf{u}} \preceq {\bf{u}}^U
\end{array},
\end{equation}

\noindent with $\bar {\bf{u}}^*_{k+1}$ then used in place of ${\bf{u}}^*_{k+1}$ in (\ref{eq:inputfilter}) and (\ref{eq:gainupper}).

If (\ref{eq:feasproj}) has a solution, then it is clear that Condition (\ref{eq:suff1}) will be satisfied by the resulting $\bar {\bf{u}}^*_{k+1}$. We now prove the feasibility of (\ref{eq:feasproj}) for the majority of practical cases.

\begin{thm}{\bf (Feasibility of Projection (\ref{eq:feasproj}))}

Consider the Jacobian matrix, ${\bf{J}}_k \in \mathbb{R}^{n_u \times n_J}$, of a simplified constraint set of (\ref{eq:feasproj}) that only consists of the gradients of the $\epsilon$-active constraints and those box constraints that are active at iteration $k$ ($i$ : $u_{k,i} = u^L_i$ or $u_{k,i} = u^U_i$). Using $(\tilde \cdot)$ to designate the inactive box constraints allows writing the constraint set of (\ref{eq:feasproj}) as\footnote{The (abuse of) notation $\tilde {\bf{u}}^L \preceq {\bf{u}} \preceq \tilde {\bf{u}}^U$ may be taken to mean generating the individual constraints of $ {\bf{u}}^L \preceq {\bf{u}} \preceq  {\bf{u}}^U$ and then removing those that are active at ${\bf u}_k$.}:

\begin{equation}\label{eq:feasre}
\begin{array}{l}
{\bf{J}}_k^T ({\bf{u}} - {\bf{u}}_{k}) \preceq \left[ {\begin{array}{*{20}c}
   -\boldsymbol{\delta}_g  \\
   {\bf{0}} \\
\end{array}} \right] \\
\tilde {\bf{u}}^L \preceq {\bf{u}} \preceq \tilde {\bf{u}}^U
\end{array}.
\end{equation}

\noindent Let ${\bf{J}}_{k,1},...,{\bf{J}}_{k,n_J}$ represent the rows of ${\bf{J}}_k^T$.

If no rows of ${\bf{J}}_k$ are negatively spanned by other rows of ${\bf{J}}_k$:

\begin{equation}\label{eq:negspan}
\nexists a_2,...,a_{n_J} \geq 0 : {\bf{J}}_{k,1} = -a_2 {\bf{J}}_{k,2} - ... - a_{n_J} {\bf{J}}_{k,{n_J}},
\end{equation}

\noindent then there exists a vector $\boldsymbol{\delta}_g \succ {\bf{0}}$ for which Projection (\ref{eq:feasproj}) has a solution.

\end{thm}

\begin{pf}

We start by noting that the proof is trivial if there are no $\epsilon$-active constraints (i.e. any point in $\mathcal{I}$ will be a feasible point).

For the general case with multiple $\epsilon$-active constraints, we use Theorem 1 to state that the top set of inequalities in (\ref{eq:feasre}) is infeasible iff:

\begin{equation}\label{eq:alt1}
\exists \boldsymbol{\nu} \in \mathbb{R}^{n_J}_+ : \hspace{1mm}\displaystyle \sum\limits_{i = 1}^{n_J} {\nu_{i}{\bf{J}}_{k,i}  } = {\bf{0}}, \hspace{3mm} \boldsymbol{\nu}^T[-\boldsymbol{\delta}_g \hspace{1mm }{\bf{0}}] < 0,
\end{equation} 

\noindent which then implies:

\begin{equation}\label{eq:lindep2}
{\bf{J}}_{k,1} = -\frac{\nu_2}{\nu_1} {\bf{J}}_{k,2} - ... - \frac{\nu_{n_J}}{\nu_1} {\bf{J}}_{k,{n_J}},
\end{equation} 

\noindent i.e. that some row is negatively spanned by the others. It follows that the top set of inequalities of (\ref{eq:feasre}) is feasible unless there is a negative spanning.

Let $\bar {\bf{u}}^*_{k+1}$ represent a solution to this set, with:

\begin{equation}\label{eq:feasre2}
{\bf{J}}_k^T(\bar {\bf{u}}^*_{k+1} - {\bf{u}}_{k}) \preceq \left[ {\begin{array}{*{20}c}
   -\boldsymbol{\delta}_g  \\
   {\bf{0}} \\
\end{array}} \right].
\end{equation}

We must now show that:

\begin{equation}\label{eq:feasre3}
\exists \bar {\bf{u}}^*_{k+1} : \tilde {\bf{u}}^L \preceq \bar {\bf{u}}^*_{k+1} \preceq \tilde {\bf{u}}^U,
\end{equation}

\noindent or that some solution of (\ref{eq:feasre2}) will hold for the inactive constraints as well. This may be restated as:

\begin{equation}\label{eq:feasre4}
\tilde {\bf{u}}^L \preceq {\bf{u}}_k + (\bar {\bf{u}}^*_{k+1} - {\bf{u}}_k) \preceq \tilde {\bf{u}}^U.
\end{equation}

Note, however, that scaling the right-hand side of (\ref{eq:feasre2}) by a scalar $\alpha > 0$ allows for the constraints to be satisfied:

\begin{equation}\label{eq:feasre5}
{\bf{J}}_k^T \alpha (\bar {\bf{u}}^*_{k+1} - {\bf{u}}_{k}) \preceq \left[ {\begin{array}{*{20}c}
   -\alpha \boldsymbol{\delta}_g  \\
   {\bf{0}} \\
\end{array}} \right],
\end{equation}

\noindent for an arbitrarily small $\alpha (\bar {\bf{u}}^*_{k+1} - {\bf{u}}_{k})$. This guarantees that a sufficiently small choice of $\alpha$ will allow:

\begin{equation}\label{eq:feasre6}
\tilde {\bf{u}}^L \preceq {\bf{u}}_k + \alpha (\bar {\bf{u}}^*_{k+1} - {\bf{u}}_k) \preceq \tilde {\bf{u}}^U,
\end{equation}

\noindent since $\tilde {\bf{u}}^L \prec {\bf{u}}_k \prec \tilde {\bf{u}}^U$, thereby completing the proof. \qed

\end{pf}

All in all, this particular result should be intuitive, as it basically excludes the case where a descent direction for some of the constraints is an ascent direction for another. As the most basic example of a case where Projection (\ref{eq:feasproj}) would be infeasible, consider two constraints, $g_{p,1}$ and $g_{p,2}$, with $\nabla g_{p,1}({\bf{u}}_k) = -\nabla g_{p,2}({\bf{u}}_k)$. This is a negative spanning and indeed matches the expected result -- one cannot go in a local descent direction for one without increasing the other.

\begin{figure}
\begin{center}
\includegraphics[width=6cm]{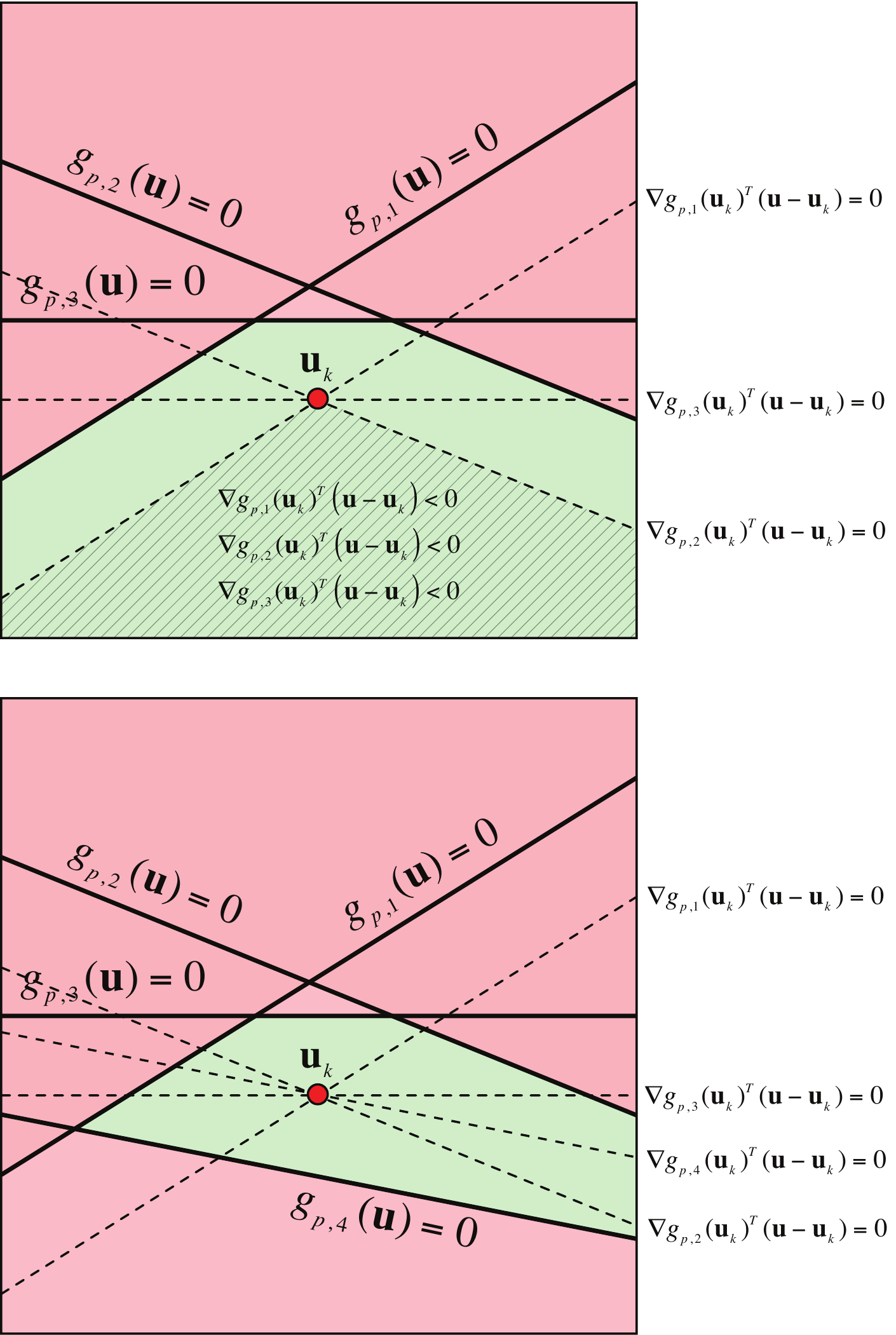}    % The printed column width  is 8.4 cm.
\caption{An illustration of the cases where (\ref{eq:feasproj}) is feasible (top) and infeasible (bottom) for a linear-constraint case. The lined region represents the space of descent directions for all of the constraints, but is empty in the bottom case due to a negative spanning between $g_{p,4}$ and the other three constraints.}
\label{fig:redcon}
\end{center}
\end{figure}

Figure {\ref{fig:redcon}} serves to illustrate this point further. In the top example, we have a simple linear case with three constraints. As their descent directions do not negatively span one another (we see, in fact, that $\nabla g_{p,1}({\bf{u}}_k)$ and $\nabla g_{p,2}({\bf{u}}_k)$ span $\nabla g_{p,3}({\bf{u}}_k)$ \emph{positively}), the projection problem is feasible (with the descent cone of all feasible projections denoted by the lined region). In the bottom case, an additional constraint, $g_{p,4}$, is added, whose descent direction lies opposite to that of the original three constraints and is spanned negatively by them. As should be clear from the figure, this renders the descent cone empty -- there is no direction that will allow for all four constraints to be decreased simultaneously. The projection is therefore infeasible in this case.

The cases where a negative spanning occurs are believed to be of little practical concern, as they represent those cases where either: (a) the choice of ${\boldsymbol \epsilon}$ is too large, or (b) the RTO problem is ill-posed to begin with. Both can be illustrated by referring to the bottom case in Figure \ref{fig:redcon}. For (a), we may assume that ${\bf{u}}_k$ is a point that is well within the feasible set, but, by an overly large choice of ${\boldsymbol \epsilon}$, we have told the projection to locally decrease all of the constraints at once. As the purpose of the projection is to avoid coming too close to the constraints that are becoming active, it is clear that such a choice is poor and unnecessary, as none of the constraints are close to active here. In this case, ${\boldsymbol \epsilon}$ should be reduced to make (\ref{eq:feasproj}) feasible. For (b), suppose that ${\boldsymbol \epsilon} \approx {\bf 0}$. If so, then the feasible space in Figure \ref{fig:redcon} is actually very, very small, with ${\bf{u}}_k$ being squeezed tightly between all four constraints, all of which are almost active -- another way to visualize this is to imagine the bottom figure in Figure \ref{fig:redcon} as having a very large zoom factor. In this case, we are probably not interested in optimization anyway, since, as would be the case with satisfying unknown equality constraints, just remaining feasible is difficult enough. For these reasons, we state that Condition (\ref{eq:suff1}) may be satisfied for the vast majority of practical cases, due to (\ref{eq:feasproj}) always having a solution in these cases provided a sufficiently low choice of ${\boldsymbol \epsilon}$ and $\boldsymbol{\delta}_g$.

The choice of projection parameters ${\boldsymbol \epsilon}$ and $\boldsymbol{\delta}_g$ deserves a few words. Perhaps the best way to present their respective roles would be by pointing out that the choice of ${\boldsymbol \epsilon}$ acts to identify the $\epsilon$-active set (and thus the local descent cone $\{{\bf{u}} : {\bf{J}}_k^T ({\bf{u}} - {\bf{u}}_{k}) \prec {\bf{0}} \}$), while $\boldsymbol{\delta}_g$ decides just how deep into this cone one would like to step.

Choosing an ${\boldsymbol \epsilon}$ that is too large can have the following negative consequences:

\begin{itemize}
\item $\{{\bf{u}} : {\bf{J}}_k^T ({\bf{u}} - {\bf{u}}_{k}) \prec {\bf{0}} \} = \varnothing$, due to a negative spanning in ${\bf{J}}_k^T$.
\item Enforcing the algorithm to stay too far away from the constraints (a local descent direction being enforced long before the constraint starts to become active).
\end{itemize}

On the contrary, choosing ${\boldsymbol \epsilon}$ to be too small will essentially cancel out the effect of Condition (\ref{eq:suff1}), as this condition will not take effect until the corresponding constraint is almost 0 (and $K_k$ almost 0). Although $K_k$ will still avoid approaching 0 asymptotically, it may become so small that the practical difference between its actual positive value and 0 will be negligible. In looser terms, this would mean ``waiting too long'' before enforcing a decrease in a given constraint value. The result in (\ref{eq:Klimit}) reflects this analytically, where the lower bound on $K_k$ clearly goes to 0 as ${\boldsymbol \epsilon} \rightarrow {\bf 0}$.

For $\boldsymbol{\delta}_g$, we remark that an arbitrarily large choice has no effect on the existence of a solution for ${\bf{J}}_k^T ({\bf{u}} - {\bf{u}}_{k}) \preceq -\boldsymbol{\delta}_g$ alone, as this is linked purely to the non-negative spanning condition of the rows of ${\bf{J}}_k^T$. However, as demonstrated in the proof of Theorem 4, choosing $\boldsymbol{\delta}_g$ to be too large would be equivalent to scaling with a very large $\alpha$, which may make satisfaction of the inactive box constraints impossible. For a small choice it is clear that the feasibility of (\ref{eq:feasproj}) will be preserved, but, by contrast, performance issues may arise. As with ${\boldsymbol \epsilon}$, we see, from a quick analytical examination of (\ref{eq:Klimit}), that the lower bound on $K_k$ goes to 0 as ${\boldsymbol \delta}_g \rightarrow {\bf 0}$, thereby allowing for very small steps and very slow progress.

\begin{figure}
\begin{center}
\includegraphics[width=12cm]{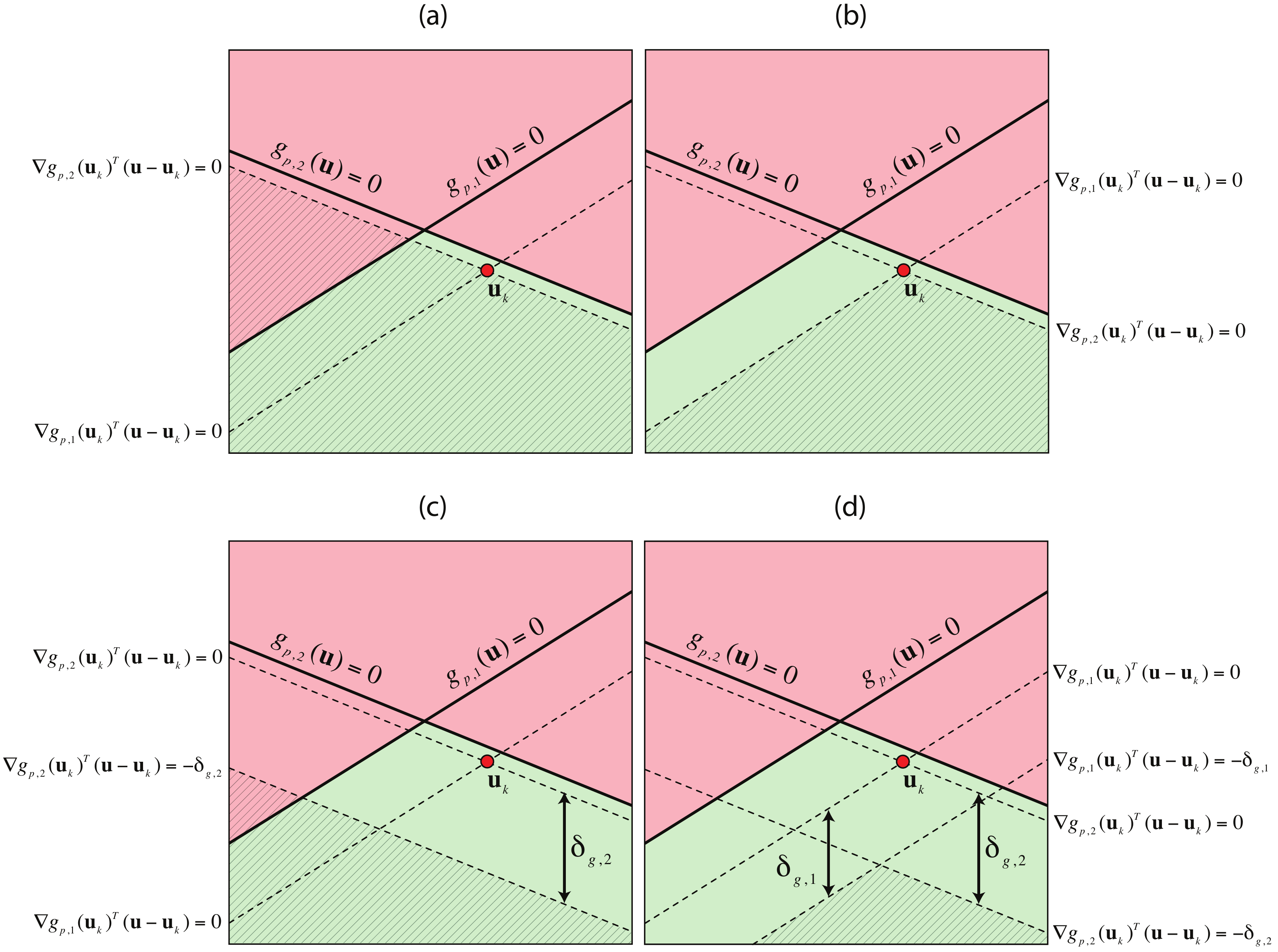}    % The printed column width  is 8.4 cm.
\caption{An illustration of how the different choices of ${\boldsymbol \epsilon}$ and $\boldsymbol{\delta}_g$ affect the feasible space of Projection (\ref{eq:feasproj}) (denoted by the lined region). (a) Small ${\boldsymbol \epsilon}$ and small $\boldsymbol{\delta}_g$ -- only $g_{p,2}$ is deemed $\epsilon$-active and only its local descent direction is enforced, (b) large ${\boldsymbol \epsilon}$ and small $\boldsymbol{\delta}_g$ -- both $g_{p,1}$ and $g_{p,2}$ are deemed $\epsilon$-active, (c) small ${\boldsymbol \epsilon}$ and large $\boldsymbol{\delta}_g$ -- only $g_{p,2}$ is deemed $\epsilon$-active and its local descent direction is enforced with a shift of $\delta_{g,2}$, (d) large ${\boldsymbol \epsilon}$ and large $\boldsymbol{\delta}_g$ -- both $g_{p,1}$ and $g_{p,2}$ are deemed $\epsilon$-active and their descent directions are enforced with shifts of $\delta_{g,1}$ and $\delta_{g,2}$. Note that the shifts ${\boldsymbol \delta}_g$ are plotted conceptually -- they are not, in general, the exact distances by which a halfspace is shifted vertically/horizontally.}
\label{fig:deltaeps}
\end{center}
\end{figure}

A geometric illustration of the different choices of ${\boldsymbol \epsilon}$ and $\boldsymbol{\delta}_g$ for a linear-constraint case is given in Figure \ref{fig:deltaeps}. We conclude by proposing that a moderate choice between the high and low extremes be used for both. However, as will be shown in Section 4, such a choice may be automated once optimality guarantees are included in the algorithm.

\subsection{Illustrative Examples}

The following RTO problem is considered:

\begin{equation}\label{eq:ex2prob}
\begin{array}{l}
\mathop {{\rm{maximize}}}\limits_{u_1,u_2} \hspace{3mm} u_2 \\
{\rm{subject}}\hspace{1mm}{\rm{to}}\hspace{3mm}g_{p,1}({\bf{u}})= u^2_1 - 0.5u_1 + u_2 -0.7 \le 0 \vspace{1mm} \\
\hspace{18mm}g_{p,2}({\bf{u}})= 2u^2_1 + 0.5u_1 + u_2 -0.75 \le 0 \vspace{1mm} \\
\hspace{18mm}u_1 \in [-0.5, 0.5], u_2 \in [0, 0.8]
\end{array},
\end{equation}

\noindent with an initial, feasible starting point of ${\bf{u}}_0 = [-0.4, 0.1]$. Here, we assume exact knowledge of the nonstrict Lipschitz matrix (the component-wise maximum of the absolute values in the Jacobian), $\tilde {\bf{K}}$, over the relevant space $\mathcal{I}$:

\begin{equation}\label{eq:lipmat}
\tilde {\bf{K}} = \left[ {\begin{array}{*{20}c}
   1.5 & 1 \\
   2.5 & 1 \\
\end{array}} \right],
\end{equation}

\noindent which we then multiply by 1.1 to obtain the strict Lipschitz constants.
 
Three algorithms are applied:

\vspace{2mm}
\noindent {\bf{Algorithm 1} -- Fixed Target}
\vspace{2mm}

By far the simplest, this ``adaptation'' involves providing ${\bf{u}}^*_{k+1}$ only once by selecting the target value of $[-0.2, 0.7]$ and keeping it there. Practically, this may correspond to a case where the initial input ${\bf{u}}_0$ has been applied for some time and has been deemed unsatisfactory, thus prompting a careful off-line analysis to choose a better operating point. However, since we do not trust to apply such an input change immediately, we would like to go there in steps, making sure that feasibility is retained throughout.

\vspace{2mm}
\noindent {\bf{Algorithm 2} -- Projected Gradient Descent with Diminishing Step}
\vspace{2mm}

Here, we adapt by taking steps along the gradient of the objective function:

\begin{equation}\label{eq:algo2}
{\bf{u}}_{k+1}^* =  {\bf{u}}_k - \frac{1}{k} \nabla \phi_p ({\bf{u}}_k),
\end{equation}

\noindent where $\nabla \phi_p ({\bf{u}}_k) = [0 \; -1]^T$. A separate projection is applied following this step for the cases when ${\bf{u}}_{k+1}^*$ falls outside of $\mathcal{I}$.

\vspace{2mm}
\noindent {\bf{Algorithm 3} -- Initial Linear Model with Constraint Adaptation}
\vspace{2mm}

A bit more involved, this algorithm employs a linear model of the constraints that is identified once in the vicinity of the initial point:

\begin{equation}\label{eq:linmod}
\begin{array}{l}
g_{1}({\bf{u}})= -1.3u_1 + u_2 -1.02 \vspace{1mm} \\
g_{2}({\bf{u}})= -1.1u_1 + u_2 -1.39
\end{array}.
\end{equation}

At each iteration, bias terms $\varepsilon$ are defined as:

\begin{equation}\label{eq:bias}
\varepsilon_{k,j} = g_{p,j}({\bf{u}}_k) - g_{j}({\bf{u}}_k),
\end{equation}

\noindent and the adaptation is then performed by solving the model-based optimization with the bias correction terms added \citep{Chachuat:08}:

\begin{equation}\label{eq:conadapt}
\begin{array}{l}
{\bf{u}}_{k+1}^* = {\rm{arg}} \mathop {{\rm{maximize}}}\limits_{u_1,u_2} \hspace{3mm} u_2 \\
\hspace{18mm}{\rm{subject}}\hspace{1mm}{\rm{to}}\hspace{3mm}g_{1}({\bf{u}}) + \varepsilon_{k,1} \le 0 \vspace{1mm} \\
\hspace{36mm}g_{2}({\bf{u}}) + \varepsilon_{k,2} \le 0 \vspace{1mm} \\
\hspace{36mm}u_1 \in [-0.5, 0.5], u_2 \in [0, 0.8]
\end{array}.
\end{equation}

\vspace{4mm}

As nominal settings, we choose $\epsilon_1 = \epsilon_2 = 0.02$ and $\delta_{g,1} = \delta_{g,2} = 0.1$ for the projection.

\begin{figure}
\begin{center}
\includegraphics[width=9.5cm]{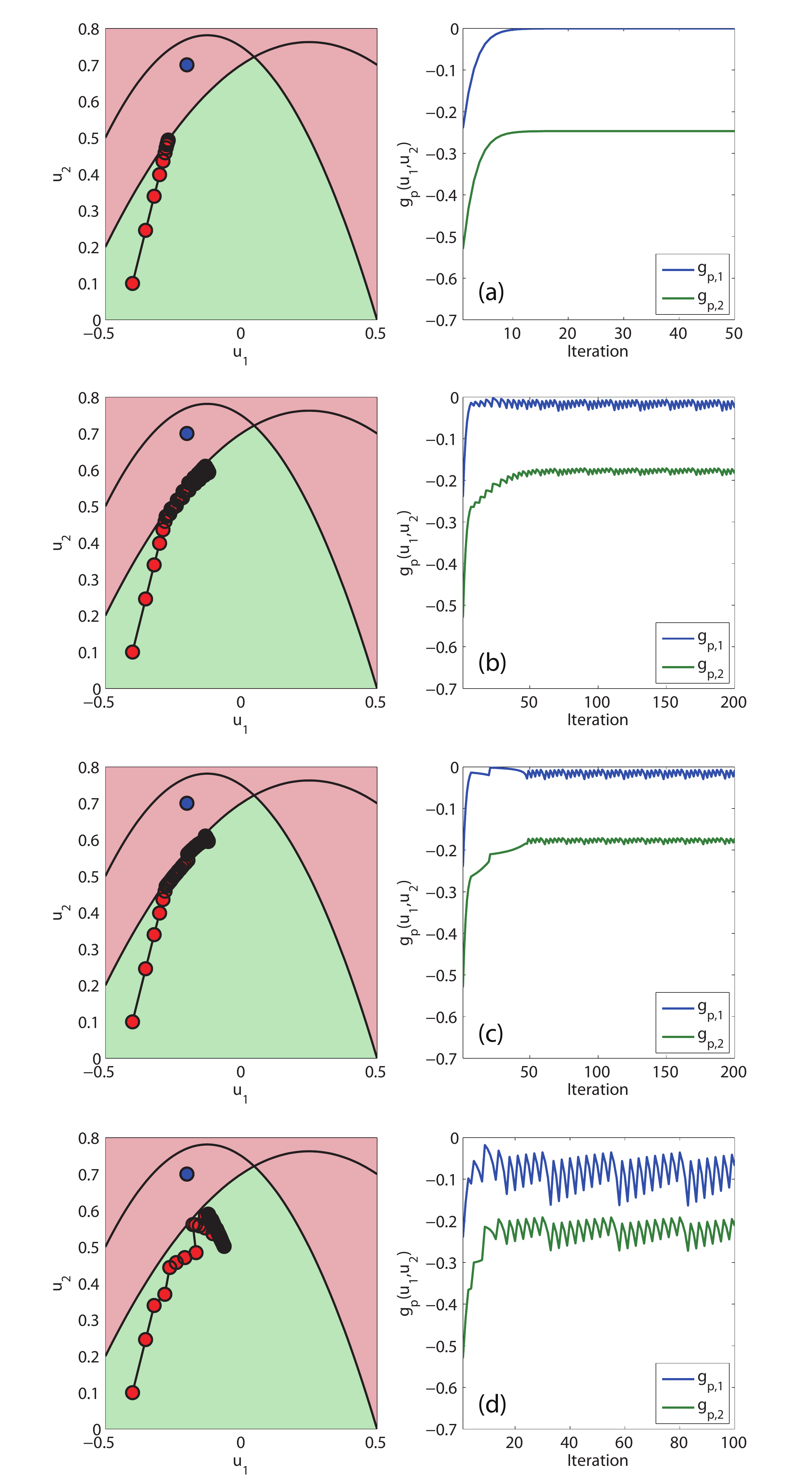}    % The printed column width  is 8.4 cm.
\caption{Feasible operation of Algorithm 1 with (a) no projection done, (b) projection with nominal ${\boldsymbol \epsilon}$ and $\boldsymbol{\delta}_{g}$ values, (c) projection with nominal ${\boldsymbol \epsilon}$ and $\delta_{g,1} = \delta_{g,2} = 0.01$, and (d) projection with $\epsilon_1 = \epsilon_2 = 0.1$ and nominal $\boldsymbol{\delta}_g$. The blue dot represents the target optimum.}
\label{fig:ex2A}
\end{center}
\end{figure}

We consider Algorithm 1 first, presenting results for the cases where no projection is done to keep a constraint from becoming active, and where this projection is applied using both the nominal and perturbed ${\boldsymbol \epsilon}$ and $\boldsymbol{\delta}_g$ values (Figure \ref{fig:ex2A}). Not enforcing Condition (\ref{eq:suff1}) leads to premature convergence, as the iterates go in a straight line towards the target point until the constraint is reached. When the projection is applied, the iterates approach the constraint but are then diverted to keep the constraint from becoming active, and the algorithm is able to reach (approximately) the closest feasible point to the target. In some sense, applying (\ref{eq:suff1}) allows us to reach a projection of the target onto the plant feasible space.

We see that varying the parameters ${\boldsymbol \epsilon}$ and $\boldsymbol{\delta}_g$ corresponds to the expectations outlined in the previous subsection. Using smaller $\boldsymbol{\delta}_g$ values clearly leads to weaker local descent -- it is seen from comparing cases (b) and (c) that the backing off from an $\epsilon$-active constraint is not as aggressive. Case (d) illustrates that choosing a larger ${\boldsymbol \epsilon}$ leads to the constraint being backed off from earlier, with the benefit of allowing larger iterate steps due to greater filter gains $K_k$. However, we see that such a choice also leads to a persistently large back-off that does not vanish.

\begin{figure}
\begin{center}
\includegraphics[width=9.5cm]{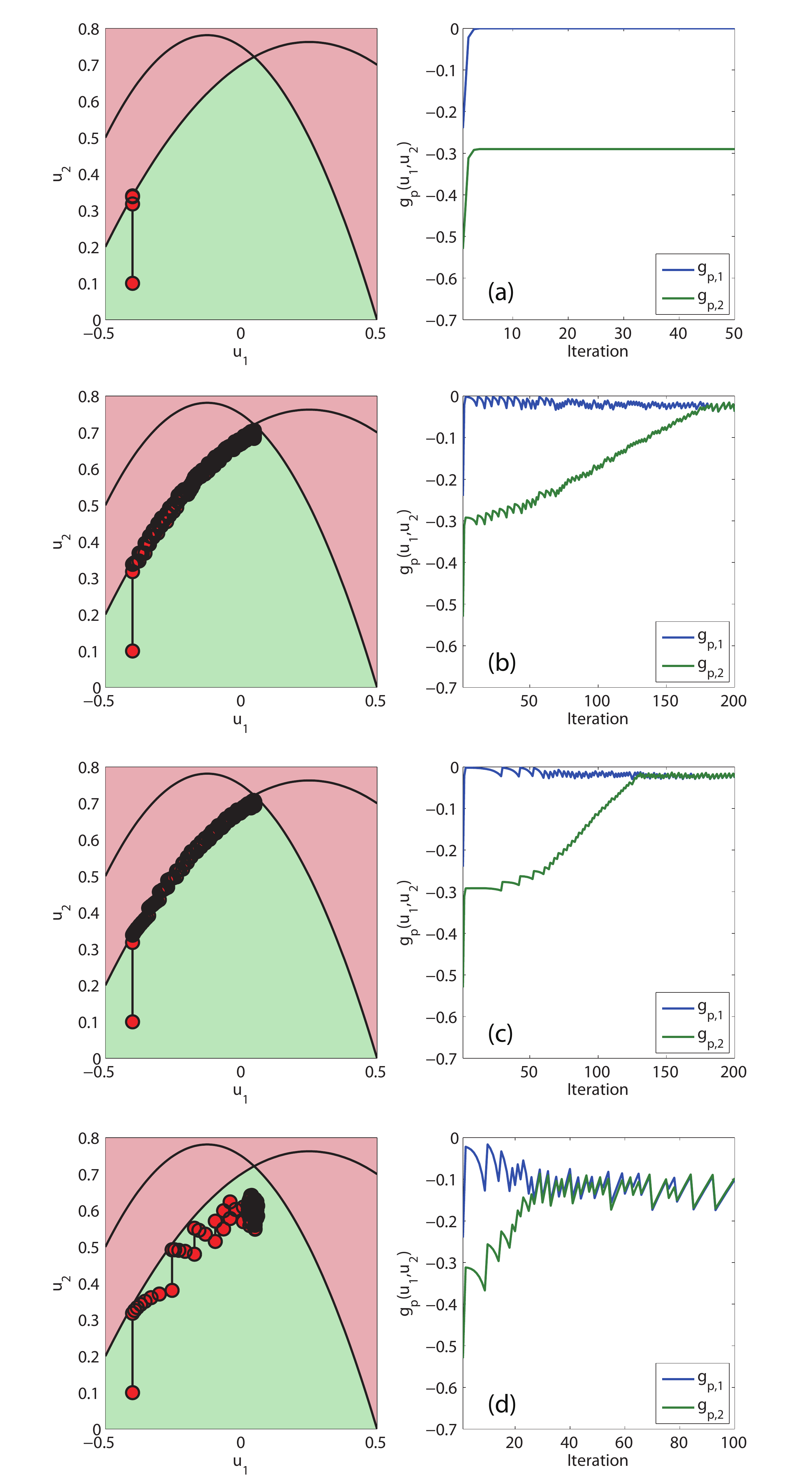}    % The printed column width  is 8.4 cm.
\caption{Feasible operation of Algorithm 2.}
\label{fig:ex2B}
\end{center}
\end{figure}

\begin{figure}
\begin{center}
\includegraphics[width=9.5cm]{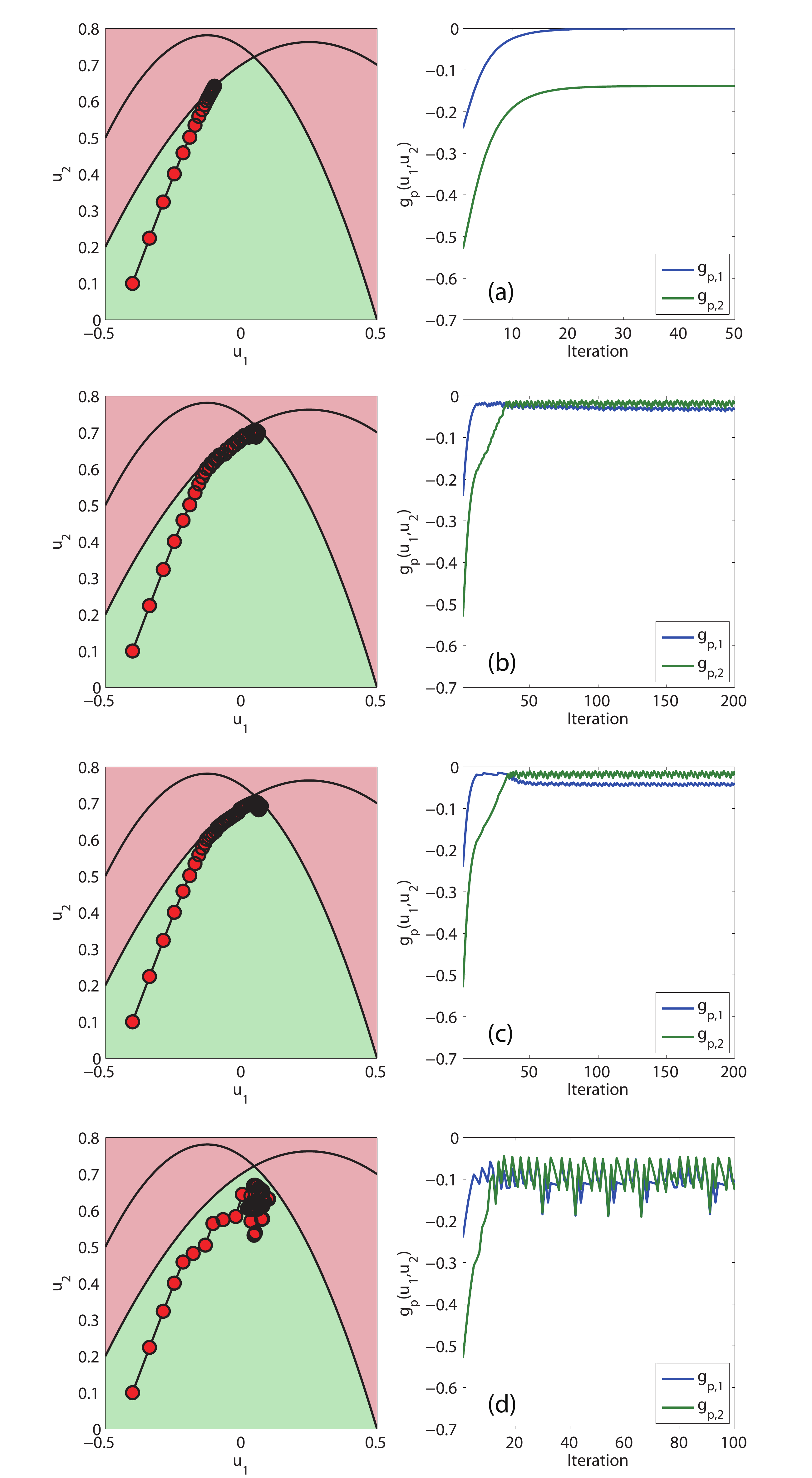}    % The printed column width  is 8.4 cm.
\caption{Feasible operation of Algorithm 3.}
\label{fig:ex2C}
\end{center}
\end{figure}

The same cases are presented (in the same order) for Algorithms 2 and 3 in Figures \ref{fig:ex2B} and \ref{fig:ex2C}, respectively. The same patterns as were noted for Algorithm 1 are seen here as well -- not applying the projection leads to convergence as soon as a constraint is reached, while applying it leads to continuous movement between iterations.

\subsection{Summary}

So as to guarantee hard constraint satisfaction of an RTO scheme at all iterations, we have reviewed some of the standard approaches to satisfy uncertain constraints. As many of these are either difficult to apply or suffer from a lack of rigor, we chose to extend the theory of the adaptive input filter method \citep{Bunin:11b} as we believe this to be a more realizable alternative. We then proceeded to highlight the main algorithmic drawback of this approach, in that it leads to premature convergence whenever one constraint becomes active and approaches 0. So as to remedy this, we proposed a condition that would be sufficient to keep the algorithm moving by projecting the target optimum onto the strict local descent halfspace of any constraint that is approaching activity. We then provided an algorithmic way of guaranteeing this condition via a projection, realized by solving a quadratic programming problem, which we then showed to be feasible for a sufficiently small choice of projection parameters ${\boldsymbol \epsilon}$ and $\boldsymbol{\delta}_g$ for the vast majority of practical cases. The effects of varying these parameters were hypothesized and corroborated in simulation, where several RTO algorithms were tried and showed the importance of having the additional condition so as to avoid premature convergence.

In the next section, we go through a similar analysis for the cost function of the RTO problem, thereby allowing us to derive conditions that, when coupled with the feasibility conditions, form the full SCFO and may be used to enforce convergence to a KKT point of the RTO problem via strictly feasible iterates.

\section{Optimality Guarantees and the Full SCFO}

We start this section by first defining ``optimality'' and then reviewing what is available in terms of optimality guarantees for the different RTO algorithms in the literature. Noting that these guarantees are algorithm-specific and rarely account for feasibility guarantees (Section 3), we propose general sufficient conditions that are algorithm-independent and preserve feasibility. As was done in the previous section, we present the mathematical analysis first, stating the sufficient conditions and proving how their presence guarantees optimality, and then provide an algorithm for their enforcement. Conceptually, the approach is very similar to the one taken in the previous section -- we simply enforce that the target optimum lie in the strict local descent halfspace of the objective function, and that the steps taken in this direction are small enough so as not to lose the descent (using, again, Lemma 3). We then finish by presenting a few examples to demonstrate the effectiveness of the theory.

\subsection{KKT Convergence and Existing Guarantees}

By ``optimality'', we simply mean KKT convergence, in that the algorithm always converges to a KKT point, ${\bf{u}}^*$, satisfying the necessary, first-order optimality conditions:

\begin{equation}\label{eq:KKT}
\begin{array}{l}
{\bf{G}}_p ({\bf{u}}^*) \preceq {\bf{0}}, \; {\bf u}^L \preceq {\bf u}^* \preceq {\bf u}^U \\
\mu_j g_{p,j} ({\bf{u}}^*) = 0, \; \zeta_i^L(u^L_i - u^*_i) = 0,\; \zeta_i^U(u^*_i - u^U_i) = 0 \vspace{1mm}\\
 \forall j = 1,...,n_g, \forall i = 1,...,n_u \\
\nabla \mathcal{L}({\bf u}^*) = \nabla \phi_p ({\bf{u}}^*) + \displaystyle \sum\limits_{j = 1}^{n_g} {\mu_j \nabla g_{p,j} ({\bf{u}}^*)} - \boldsymbol{\zeta}^L + \boldsymbol{\zeta}^U = {\bf{0}}
\end{array},
\end{equation}

\noindent with $\boldsymbol{\mu} \in \mathbb{R}^{n_g}_+$ and $\boldsymbol{\zeta}^L, \boldsymbol{\zeta}^U \in \mathbb{R}^{n_u}_+$ vectors of Lagrange multipliers corresponding to the uncertain and box constraints, respectively, and $\mathcal{L}({\bf u}) : \mathbb{R}^{n_u} \rightarrow \mathbb{R}^{n_u}$ the Lagrangian function that defines the KKT stationarity condition.

While this is not \emph{sufficient} to guarantee that the KKT point is a local minimum (a second-order condition is also required), we limit ourselves to this definition as the KKT points that are not local minima (e.g. local maxima and saddle points) will not be algorithmically stable in any descent method. As the conditions proposed in this work enforce a descent method, we will consider (\ref{eq:KKT}) to be ``practically sufficient'' to guarantee local minimality, even though we are unable to guarantee this rigorously -- while conceptually possible, such guarantees would require much stronger assumptions (knowledge of the plant Hessian) that may be extremely difficult to uphold in practice.

In considering different RTO algorithms, one usually sees a clear distinction between model-based and model-free methods with regard to KKT convergence -- the conditions for the latter are usually much simpler and easier to realize (though the algorithms are slower and require more iterations), while the former, despite generally faster convergence, have conditions that require difficult assumptions on the model used and the potential plant-model error. The guarantees available for the different approaches, as well as the price to pay for such guarantees, are summarized as follows:

\begin{itemize}
\item The two-step parameter-identification-and-optimization approach, which is the natural and standard approach to optimize many processes, has very weak guarantees for converging to an actual KKT point of the plant \citep{Forbes:94a,Brdys:05, Marchetti:09, Quelhas:12}. In the general case, such convergence may only be guaranteed if there exists a set of parameters for which the model has the same KKT point as the plant \citep{Biegler:85} and if such a set is identified by the parameter estimation method. When this is not the case, it is generally accepted that the converged point of the model optimization will not be a plant KKT point. In the rare case when structural uncertainty is not present and all of the uncertainty is described by the parameters, the convergence properties of the two-step approach are very simple -- it is sufficient to identify the correct parameters, to optimize the model once, and to apply the result to the plant.
\item Methods that employ measurements to correct the 0$^{th}$- and 1$^{st}$-order error of the model, known as \emph{ISOPE} \citep{Brdys:05} or \emph{modifier adaptation} \citep{Marchetti:09a}, are guaranteed to be at a KKT point \emph{upon convergence}, which is an attractive property but nevertheless requires that the scheme converge. Local necessary conditions for convergence have been detailed \citep{Marchetti:09a}, but are difficult to implement as they require the knowledge of the first and second derivatives of the plant at the KKT point. Sufficient conditions for the case without uncertain inequality constraints are also available \citep[Th. 4.1]{Brdys:05}, but depend on quantities that may not be easy to compute (i.e. the global minimal eigenvalue of the model Hessian and its relation to the quadratic upper bound of the plant) and assume convexity in the constraints that are known. Sufficient conditions for the case with both uncertain cost and constraints have been proposed \citep{Chachuat:08b}, but are largely abstract and require multiple quadratic upper bounds that are linked to the plant-model error, as well as certain continuity assumptions on the model. 
\item Model-free methods such as the simplex, direct search, and gradient descent come with simple guarantees of KKT convergence that rely almost entirely on some sort of line search \citep{Fletcher:87,Conn:09}. These are easy to implement as they need only continued experimentation and sufficiently small steps, but have the obvious drawback of inefficiency and slow convergence as they do not use \emph{a priori} knowledge. 
\end{itemize} 

It is important to note that none of these guarantees, to the authors' best knowledge, take hard constraints into account. As such, while guaranteeing convergence to a KKT point may be possible using the available theory for a specific algorithm, doing so in the presence of hard constraints is not something that has been addressed.

\subsection{General Sufficient Conditions for Feasibility and Optimality}

We begin by proposing a set of sufficient conditions for monotonic cost decrease in a general RTO algorithm.

\begin{lem}{\bf (Minimal Cost Decrease Between Iterations)}

Let the following two conditions be satisfied at every iteration of the RTO algorithm:

\begin{equation}\label{eq:suff2}
\begin{array}{l}
\nabla \phi_p({\bf{u}}_k)^T ({\bf{u}}^*_{k+1} - {\bf{u}}_{k}) \leq -\delta_\phi \vspace{2mm}\\
0 < K_{min} \leq K_k \leq \\
\hspace{15mm}-2 \displaystyle \frac{\nabla \phi_p({\bf{u}}_k)^T ({\bf{u}}^*_{k+1} - {\bf{u}}_{k})}{({\bf{u}}^*_{k+1} - {\bf{u}}_{k})^T \overline {\bf{Q}}_\phi ({\bf{u}}^*_{k+1} - {\bf{u}}_{k})} - K_{min}
\end{array},
\end{equation}

\noindent with $\delta_\phi > 0$, $\overline {\bf{Q}}_\phi \succ 0$ the quadratic upper bound of $\phi_p ({\bf{u}})$ as defined in Lemma 2, and $K_{min} > 0$ some minimal value achieved by the adaptive filter gain $K_k$.\footnote{For (\ref{eq:suff2}) to be true, it is implicit that $K_{min} \leq - \frac{\nabla \phi_p({\bf{u}}_k)^T ({\bf{u}}^*_{k+1} - {\bf{u}}_{k})}{({\bf{u}}^*_{k+1} - {\bf{u}}_{k})^T \overline {\bf{Q}}_\phi ({\bf{u}}^*_{k+1} - {\bf{u}}_{k})}$.}

Then, the change in the cost from iteration to iteration will be bounded as:

\begin{equation}\label{eq:costchange}
\begin{array}{l}
\phi_p ({\bf{u}}_{k+1}) - \phi_p ({\bf{u}}_{k}) \leq -K_{min} \; \delta_\phi + \vspace{1mm}\\
\hspace{10mm}\displaystyle \frac{K^2_{min}}{2} ({\bf{u}}^{U} - {\bf{u}}^{L})^T \overline {\bf{Q}}_\phi ({\bf{u}}^{U} - {\bf{u}}^{L})  < 0
\end{array}.
\end{equation}

\end{lem}

\begin{pf}

Using Lemma 3, we know that the worst-case evolution (\ref{eq:Qboundw2}), being strictly convex with respect to $K_k$, achieves its maximum at either $K_k = K_{min}$ or $K_k = -2\frac{\nabla \phi_p({\bf{u}}_k)^T ({\bf{u}}^*_{k+1} - {\bf{u}}_{k})}{({\bf{u}}^*_{k+1} - {\bf{u}}_{k})^T \overline {\bf{Q}}_\phi ({\bf{u}}^*_{k+1} - {\bf{u}}_{k})}$ $- K_{min}$ and has a value strictly less than 0.

Substituting either of these bounds into the worst-case evolution expression (\ref{eq:Qboundw2}) for $\phi_p ({\bf{u}}_{k+1}) - \phi_p ({\bf{u}}_{k})$ yields:

\begin{equation}\label{eq:proofD3}
\begin{array}{l}
\phi_p({\bf{u}}_{k+1}) - \phi_p({\bf{u}}_k) \leq  K_{min} \nabla \phi_p({\bf{u}}_k)^T ({\bf{u}}^*_{k+1} - {\bf{u}}_{k}) + \vspace{1mm}\\
\hspace{35mm}\displaystyle \frac{K^2_{min}}{2}({\bf{u}}^*_{k+1} - {\bf{u}}_{k})^T \overline {\bf{Q}}_\phi ({\bf{u}}^*_{k+1} - {\bf{u}}_{k}) 
\end{array},
\end{equation}

\noindent which, with the worst-case values of the linear and quadratic terms, leads to the bound in (\ref{eq:costchange}). \qed

\end{pf}

Conditions (\ref{eq:suff2}) are useful in that they give a guarantee of a minimal cost decrease between two consecutive iterations. As the cost is, by assumption, continuous over a compact domain and thereby bounded, it follows that meeting such conditions indefinitely is not possible. The following gives a global upper bound on the maximum number of iterations for which these conditions may be met.

\begin{corol}{\bf (Upper Bound on Number of Cost-Decreasing, Feasible Iterations)}

Let $\phi_{p,min} = {\rm{min}} (\phi_p({\bf{u}}) : {\bf{u}} \in \mathcal{I},  {\bf{G}}_p({\bf{u}}) \preceq {\bf{0}} )$ and let the conditions of Lemma 4 hold whenever possible for some fixed $\delta_\phi > 0$, with all iterates satisfying the hard inequality and box constraints. The number of iterations for which the conditions of Lemma 4 are satisfied cannot exceed:

\begin{equation}\label{eq:upperiter}
\frac{\phi_p({\bf{u}}_0)-\phi_{p,min}}{K_{min} \; \delta_\phi - \displaystyle \frac{K^2_{min}}{2} ({\bf{u}}^{U} - {\bf{u}}^{L})^T \overline {\bf{Q}}_\phi ({\bf{u}}^{U} - {\bf{u}}^{L}) } < \infty.
\end{equation}

\end{corol}

\begin{pf}

This results directly from Lemma 4. As $\phi_p({\bf{u}}_0)-\phi_{p,min}$ represents the greatest possible suboptimality gap between the initial point and the global minimum, the validity of Conditions (\ref{eq:suff2}) for more than this number of iterations would reduce the cost below its global minimum, which is not possible without losing feasibility. \qed

\end{pf}

We now combine all of the results obtained so far to give the full SCFO.

\begin{thm}{\bf (Sufficient Conditions for Feasible-Side Convergence to a KKT Point)}

Let an RTO algorithm satisfy Conditions (\ref{eq:gainupper}), (\ref{eq:suff1}), and (\ref{eq:suff2}) for all iterations when it is possible to do so for some fixed choice of ${\boldsymbol \epsilon}$, $\boldsymbol{\delta}_g$, $\delta_\phi \succ {\bf{0}}$, and yield ${\bf u}_{k+1} = {\bf u}_k$ otherwise. Defining the KKT error $\mathcal{E}$ as the minimal sum of squared errors in the stationarity and complementary slackness conditions of (\ref{eq:KKT}):

\begin{equation}\label{eq:KKTerr}
\begin{array}{l}
\mathcal{E}({\bf u}) = \mathop {\inf} \limits_{{\boldsymbol \mu}, {\boldsymbol \zeta^L}, {\boldsymbol \zeta^U} \succeq {\bf 0}} \Big ( \nabla \mathcal{L}({\bf u})^T \nabla \mathcal{L}({\bf u}) +  \displaystyle \sum\limits_{j = 1}^{n_g} \left[ \mu_j g_{p,j} ({\bf u}) \right]^2 + \\
\hspace{30mm}\displaystyle \sum\limits_{i = 1}^{n_u} \left[ ( \zeta^L_i (u^L_i - u_i) )^2 + ( \zeta^U_i (u_i - u^U_i) )^2  \right] \Big )
\end{array},
\end{equation}

\noindent it follows that such an algorithm will:

\begin{enumerate}[(i)]
\item Converge to a static point, ${\bf u}_\infty$, in a finite number of iterations.
\item Satisfy the plant constraints at every iteration.
\item Decrease the cost at every iteration until ${\bf u}_\infty$.
\item Have the KKT error at ${\bf u}_\infty$ go to 0 in the limit with respect to ${\boldsymbol \epsilon}$, $\boldsymbol{\delta}_g$, and $\delta_\phi$:

\begin{equation}\label{eq:KKTerr0}
\mathop {\lim} \limits_{{\boldsymbol \epsilon}, \boldsymbol{\delta}_g, \delta_\phi \rightarrow {\bf 0}} \mathcal{E}({\bf u}_\infty) = 0.
\end{equation}

\end{enumerate}

\end{thm}

\begin{pf}

We refer the reader to the appendix. \qed

\end{pf}

Conditions (\ref{eq:gainupper}), (\ref{eq:suff1}), and (\ref{eq:suff2}), together with ${\bf u}_{k+1}^* \in \mathcal{I}$, constitute the full SCFO -- guaranteeing finite-time convergence to a stationary point with arbitrarily small KKT error as ${\boldsymbol \epsilon}$, $\boldsymbol{\delta}_g$, and $\delta_\phi$ are made arbitrarily small. We note that we cannot decouple optimality from feasibility, as we must have the latter to guarantee the former. However, it may be shown that convergence via infeasible iterates is also possible via an additional implementation technique -- this will be addressed in the companion work \citep{Bunin:12c}.

\begin{figure}
\begin{center}
\includegraphics[width=6cm]{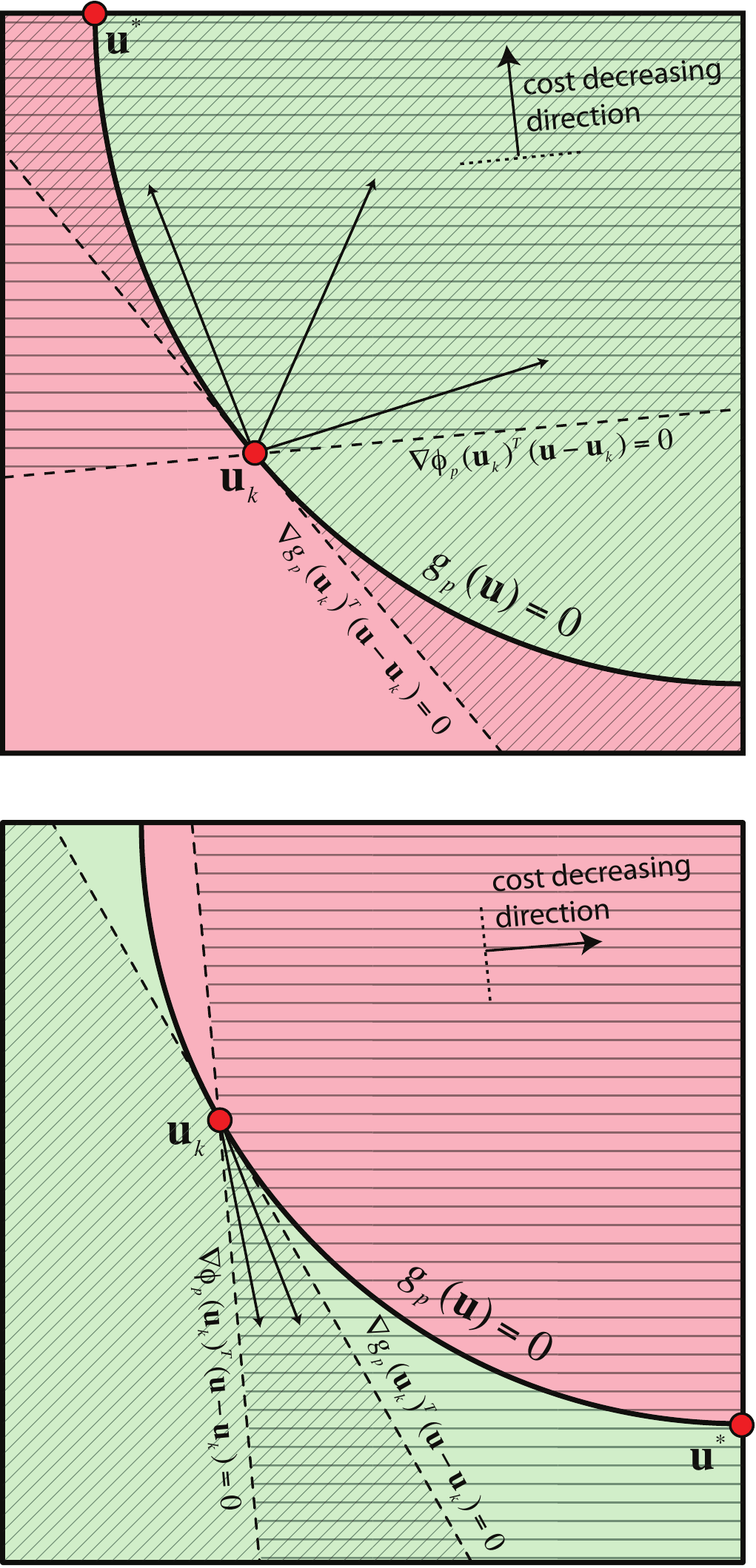}    % The printed column width  is 8.4 cm.
\caption{The geometric interpretation of the SCFO for the earlier examples of Figure \ref{fig:feasgraph2} with a linear cost included. Here, the diagonally-lined regions represent the local descent halfspaces for the constraints, while the horizontally-lined regions represent the local descent halfspaces for the cost. The double-lined regions (cones) are therefore those that satisfy both Conditions (\ref{eq:suff1}) and (\ref{eq:suff2}), with the arrows showing potential directions that may then be taken by the RTO algorithm in accordance to the SCFO. It is not difficult to visualize that such a cone will always exist unless the current iterate ${\bf u}_k$ is already a KKT point.}
\label{fig:feascongraph}
\end{center}
\end{figure}

We give a geometrical illustration of the combined effect of Conditions (\ref{eq:suff1}) and (\ref{eq:suff2}) in Figure \ref{fig:feascongraph} (as an extension of Figure \ref{fig:feasgraph2}), which allows for a much simpler interpretation -- by enforcing that the RTO algorithm always move in a locally cost-decreasing and feasible direction, it is natural that it converge when there no longer exists such a direction to move in (i.e. to the geometric definition of a KKT point).

\subsection{Basic Implementation}

We propose a similar implementation as we did for the feasibility-only case, and perform the following projection:

\begin{equation}\label{eq:feasproj2}
\begin{array}{l}
\bar {\bf{u}}^*_{k+1} = {\rm{arg}} \mathop {{\rm{minimize}}}\limits_{{\bf{u}}}\hspace{4mm}\left\| {\bf{u}} - {\bf{u}}^*_{k+1} \right\|_2^2 \vspace{1mm}  \\
\hspace{17mm}{\rm{subject}}\hspace{1mm}{\rm{to}}\hspace{3mm}\nabla g_{p,j}({\bf{u}}_k)^T ({\bf{u}} - {\bf{u}}_{k}) \leq -\delta_{g,j} \vspace{1mm} \\
\hspace{36mm}\forall j: g_{p,j}({\bf{u}}_k) \geq -\epsilon_j \vspace{1mm}\\
\hspace{36mm}\nabla \phi_{p}({\bf{u}}_k)^T ({\bf{u}} - {\bf{u}}_{k}) \leq -\delta_\phi \vspace{1mm} \\
\hspace{36mm}{\bf{u}}^L \preceq {\bf{u}} \preceq {\bf{u}}^U
\end{array},
\end{equation}

\noindent followed by applying (\ref{eq:inputfilter}) with respect to $\bar {\bf u}_{k+1}^*$, with the filter gain defined as:

\begin{equation}\label{eq:gainuppereq2}
\begin{array}{l}
K_{k} := \mathop {\min} \Bigg \{ \mathop{\min}\limits_{j = 1,...,n_g} \left[\displaystyle \frac{{-g_{p,j}({\bf{u}}_{k } )}}{\displaystyle \sum\limits_{i = 1}^{n_u} {\kappa_{ji} | \bar u^*_{k+1,i} - u_{k,i} | }} \right],  \\
\hspace{25mm} -1.99 \displaystyle \frac{\nabla \phi_p({\bf{u}}_k)^T (\bar {\bf{u}}^*_{k+1} - {\bf{u}}_{k})}{(\bar {\bf{u}}^*_{k+1} - {\bf{u}}_{k})^T \overline {\bf{Q}}_\phi (\bar {\bf{u}}^*_{k+1} - {\bf{u}}_{k})} \Bigg \} \vspace{2mm}\\
K_k > 1 \rightarrow K_k := 1
\end{array}\;\;\;\;,
\end{equation}

\noindent where we use 1.99 instead of 2 as an approximation of the inequality in (\ref{eq:suff2}).

As shown in Theorem 5, Problem (\ref{eq:feasproj2}) will have a solution for ${\boldsymbol \epsilon}$, $\boldsymbol{\delta}_g$, and $\delta_\phi$ sufficiently small unless ${\bf{u}}_k$ is already a KKT point. As might be expected, and as was already shown in both Figure \ref{fig:deltaeps} and the illustrative example in Section 3, the exact way in which the projection is done (i.e. the values of ${\boldsymbol \epsilon}$, $\boldsymbol{\delta}_g$, and $\delta_\phi$) will influence the performance of the algorithm -- we give an extension of the previous geometric illustration, adding a linear cost, in Figure \ref{fig:deltaeps2}. It is not clear what choices would lead to the best performance, but the following general trends may be expected and have been observed:

\begin{figure}
\begin{center}
\includegraphics[width=12cm]{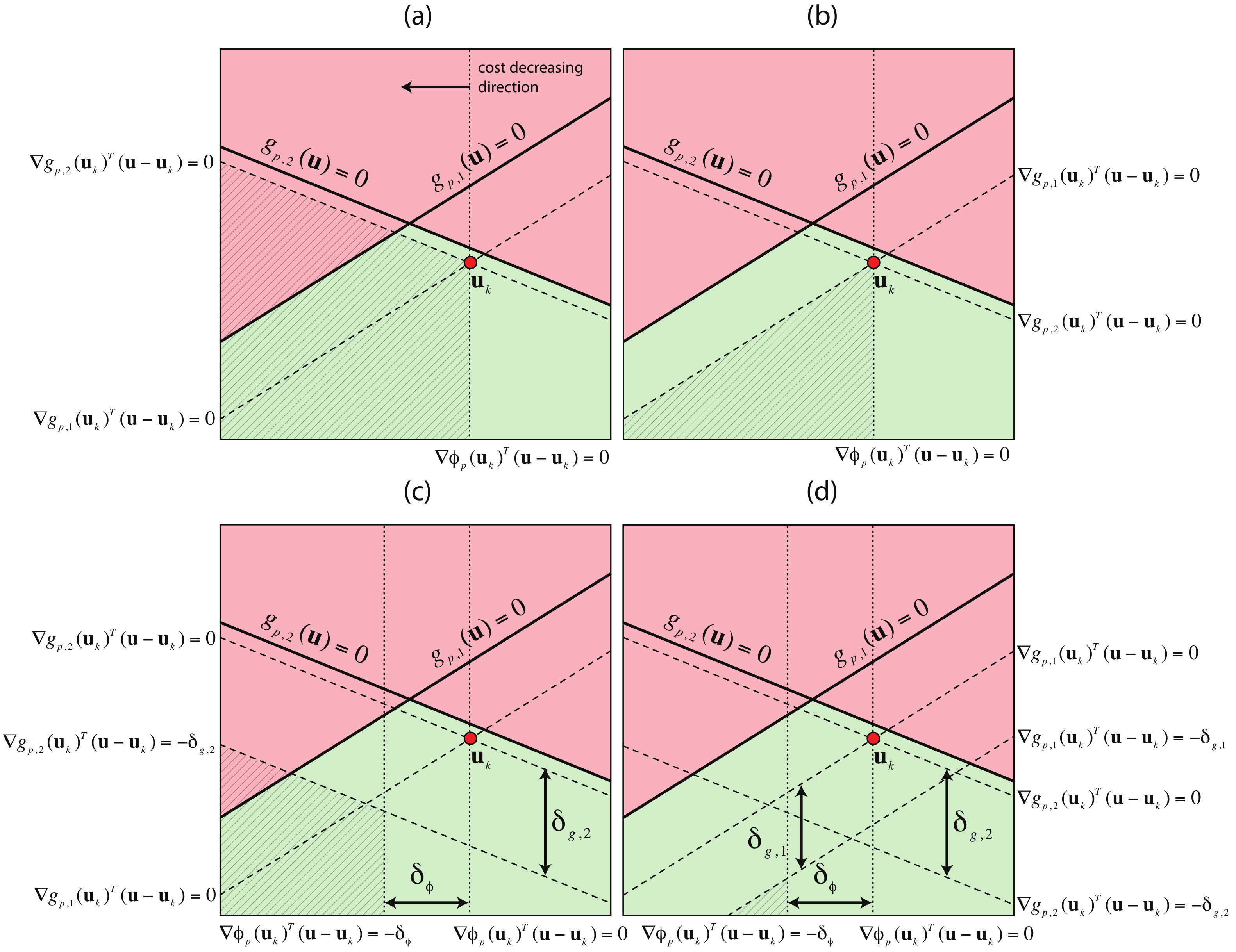}    % The printed column width  is 8.4 cm.
\caption{An illustration of how the different choices of ${\boldsymbol \epsilon}$, $\boldsymbol{\delta}_g$, and $\delta_\phi$ affect the feasible space of Projection (\ref{eq:feasproj2}) (shown via the lined regions). (a) Small ${\boldsymbol \epsilon}$ and small $\boldsymbol{\delta}_g$, $\delta_\phi$, (b) large ${\boldsymbol \epsilon}$ and small $\boldsymbol{\delta}_g$, $\delta_\phi$, (c) small ${\boldsymbol \epsilon}$ and large $\boldsymbol{\delta}_g$, $\delta_\phi$, (d) large ${\boldsymbol \epsilon}$ and large $\boldsymbol{\delta}_g$, $\delta_\phi$. Note that the shifts ${\boldsymbol \delta}_g$ and $\delta_\phi$ are plotted conceptually -- they are not, in general, the exact distances by which a halfspace is shifted vertically/horizontally.}
\label{fig:deltaeps2}
\end{center}
\end{figure}

\begin{itemize}
\item Small constant values of ${\boldsymbol \epsilon}$, $\boldsymbol{\delta}_g$, and $\delta_\phi$ approximate the true SCFO well, i.e. the $\epsilon$-active constraints are close to the true active constraints, and Conditions (\ref{eq:suff1}) and (\ref{eq:suff2}) are close to the strict local descent halfspace conditions (the shifts seen in Cases (c) and (d) of Figure \ref{fig:deltaeps2} are negligible). This leads to ``small'' KKT error upon convergence as characterized analytically in Theorem 5. However, the conditions may be enforced ``too late'', with a constraint becoming close to active, and the amount of local descent may be small. This may result in very small values for $K_k$ and therefore small steps, causing the algorithm to progress very slowly, especially when close to a constraint. The number of worst-case iterations as given by (\ref{eq:upperiter}) also grows accordingly, and the cost improvement guaranteed at each step, via Lemma 4, lessens.
\item Large constant values of ${\boldsymbol \epsilon}$, $\boldsymbol{\delta}_g$, and $\delta_\phi$ have the opposite effect. A greater $\epsilon$-active set is chosen due to larger ${\boldsymbol \epsilon}$ and thus the algorithm tries to stay away from more constraints and tries to stay away from them earlier, long before they become close to active. The amount of local descent in both the cost and $\epsilon$-active constraints is also greater, and so it should be expected that the algorithm will aim for larger steps (this is seen clearly in Case (d) of Figure \ref{fig:deltaeps2}). However, the SCFO are poorly approximated and the worst-case KKT error upon convergence may be large. Clearly, large $\boldsymbol{\epsilon}$, $\boldsymbol{\delta}_g$, and $\delta_\phi$ all increase the chance of the conditions being impossible to satisfy -- for example, were $\boldsymbol{\delta}_g$ and $\delta_\phi$ increased further in Case (d) of Figure \ref{fig:deltaeps2}, the projection would clearly become infeasible as the small triangular feasible space would simply vanish.
\end{itemize}

As such, a sound approach would be to use larger values of $\boldsymbol{\epsilon}$, $\boldsymbol{\delta}_g$, and $\delta_\phi$ if possible, as this appears to lead to faster progress, and to use smaller values otherwise to enforce true KKT convergence, which is guaranteed only as $\boldsymbol{\epsilon},\boldsymbol{\delta}_g, \delta_\phi \rightarrow {\bf 0}$. To this end, we propose the following method to auto-select and auto-tune these parameters:

\vspace{2mm}
\noindent {\bf{Initialization -- Done Prior to RTO}}
\vspace{2mm}

\begin{enumerate}
\item The degrees of the different constraints' activity, as well as the amount of local descent in the constraints and cost, should be on a comparable scale. A simple way to ensure this is by scaling the constraints and cost with respect to their ranges (e.g. if $\mathop{\min}\limits_{{\bf{u}} \in \mathcal{I}} g_p({\bf{u}}) = -100$, the scaled constraint may be re-defined as $g_p({\bf{u}}) := 0.01g_p({\bf{u}})$). 

Setting upper and lower limits on the projection parameters, let $\overline {\boldsymbol \epsilon} =  \boldsymbol{\overline \delta}_g = {\bf{1}}$, $\overline \delta_\phi = 1$, and choose $\underline {\boldsymbol \epsilon}$, $\boldsymbol{\underline \delta}_g$, and $\underline \delta_\phi$ sufficiently small (e.g. 10$^{-6}$) so that the approximation error of the active set by the $\underline \epsilon$-active set and of the strict local descent conditions by the nonstrict inequality conditions, with $-\boldsymbol{\underline \delta}_g$ and $-\underline \delta_\phi$ on the right-hand sides, is negligible. 
\end{enumerate}

\vspace{2mm}
\noindent {\bf{Search for a Feasible Projection -- Before Each RTO Iteration}}
\vspace{2mm}

\begin{enumerate}
\setcounter{enumi}{1}
\item Set ${\boldsymbol \epsilon} := \overline {\boldsymbol \epsilon}$, $\boldsymbol{\delta}_g := \boldsymbol{\overline \delta}_g$, and $\delta_\phi := \overline \delta_\phi$.
\item Check the feasibility of (\ref{eq:feasproj2}) for the given choice of ${\boldsymbol \epsilon}$, $\boldsymbol{\delta}_g$, and $\delta_\phi$ \footnote{As the arguably best way to verify the feasibility of (\ref{eq:feasproj2}) is to solve a linear programming optimization problem, we write it in this form, with the ``minimization of 0'' used to denote that no optimization actually takes place.}:

\begin{equation}\label{eq:feasproj2slack}
\begin{array}{l}
\hspace{1mm}\mathop {{\rm{minimize}}}\limits_{{\bf{u}}} \hspace{3mm} 0 \vspace{1mm}  \\
{\rm{subject}}\hspace{1mm}{\rm{to}}\hspace{3mm}\nabla g_{p,j}({\bf{u}}_k)^T ({\bf{u}} - {\bf{u}}_{k}) \leq -\delta_{g,j} \vspace{1mm} \\
\hspace{19mm}\forall j: g_{p,j}({\bf{u}}_k) \geq -\epsilon_j \vspace{1mm}\\
\hspace{19mm}\nabla \phi_{p}({\bf{u}}_k)^T ({\bf{u}} - {\bf{u}}_{k}) \leq -\delta_\phi \vspace{1mm} \\
\hspace{19mm}{\bf{u}}^L \preceq {\bf{u}} \preceq {\bf{u}}^U
\end{array}.
\end{equation}

If (\ref{eq:feasproj2slack}) does not have a solution, set ${\boldsymbol \epsilon} := 0.5 {\boldsymbol \epsilon}$, $\boldsymbol{\delta}_g := 0.5 \boldsymbol{\delta}_g$, $\delta_\phi := 0.5 \delta_\phi$, and attempt to re-solve (\ref{eq:feasproj2slack}). Otherwise, proceed to Step 4.

\item If ${\boldsymbol \epsilon} \succeq \underline {\boldsymbol \epsilon}$, $\boldsymbol{\delta}_g \succeq \boldsymbol{\underline \delta}_g$, or $\delta_\phi \geq \underline \delta_\phi$, solve (\ref{eq:feasproj2}) with ${\boldsymbol \epsilon}$, $\boldsymbol{\delta}_g$, and $\delta_\phi$. Else, terminate (Step 5).

\end{enumerate}

\vspace{2mm}
\noindent {\bf{Termination -- Declared Convergence to KKT Point}}
\vspace{2mm}

\begin{enumerate}
\setcounter{enumi}{4}
\item If ${\boldsymbol \epsilon} \prec \underline {\boldsymbol \epsilon}$, $\boldsymbol{\delta}_g \prec \boldsymbol{\underline \delta}_g$, and $\delta_\phi < \underline \delta_\phi$, terminate, with $\bar {\bf{u}}^*_{k+1} := {\bf{u}}_{k}$.
\end{enumerate}

The basic philosophy of this method is in using larger values of ${\boldsymbol \epsilon}$, $\boldsymbol{\delta}_g$, and $\delta_\phi$ to both stay away from the constraints (large ${\boldsymbol \epsilon}$ and $\boldsymbol{\delta}_g$) and to significantly decrease the cost (large $\delta_\phi$) \emph{when possible}. Once this becomes impossible due to feasibility issues in (\ref{eq:feasproj2}), these parameters are lowered until (\ref{eq:feasproj2}) becomes feasible. From Theorem 5, it is clear that infeasibility that persists even when ${\boldsymbol \epsilon} \prec \underline {\boldsymbol \epsilon}$, $\boldsymbol{\delta}_g \prec \boldsymbol{\underline \delta}_g$, and $\delta_\phi < \underline \delta_\phi$ implies that convergence to a KKT point has been achieved with very good accuracy, provided that $\underline {\boldsymbol \epsilon}$, $\underline {\boldsymbol \delta}_g$, and $\underline \delta_{\phi}$ are not too large.

We finish by remarking that the update laws of dividing ${\boldsymbol \epsilon}$, $\boldsymbol{\delta}_g$, and $\delta_\phi$ by 2 in Step 3 are purely heuristic. One could propose other, perhaps better performing, rules, but such an optimization of the reduction law for ${\boldsymbol \epsilon}$, $\boldsymbol{\delta}_g$, and $\delta_\phi$ (itself a potential RTO problem) is outside the scope of the present work. As will be seen in the next subsection, the proposed law gives satisfactory results for the cases studied here.

\subsection{Illustrative Examples}

Consider the following RTO problem with one convex and two concave uncertain constraints:

\begin{equation}\label{eq:ex4prob}
\begin{array}{l}
\mathop {{\rm{minimize}}}\limits_{u_1,u_2} \hspace{3mm} \phi_{p}({\bf{u}}) = (u_1-0.5)^2 + (u_2-0.4)^2 \\
{\rm{subject}}\hspace{1mm}{\rm{to}}\hspace{3mm}g_{p,1}({\bf{u}})= -6u^2_1 - 3.5u_1 + u_2 -0.6 \le 0 \vspace{1mm} \\
\hspace{18mm}g_{p,2}({\bf{u}})= 2u^2_1 + 0.5u_1 + u_2 -0.75 \le 0 \vspace{1mm} \\
\hspace{18mm}g_{p,3}({\bf{u}})= -u^2_1 - (u_2-0.15)^2 +0.01 \le 0 \vspace{1mm} \\
\hspace{18mm}u_1 \in [-0.5, 0.5], u_2 \in [0, 0.8]
\end{array},
\end{equation}

\noindent for which the matrix of Lipschitz constants is defined as:

\begin{equation}\label{eq:lipmat2}
{\bf{K}} = 1.1\left[ {\begin{array}{*{20}c}
   9.5 & 1 \\
   2.5 & 1 \\
   1 & 1.3 \\
\end{array}} \right],
\end{equation}

\noindent and $\overline {\bf{Q}}_\phi$ chosen as the Hessian of the (uncertain) cost function:

\begin{equation}\label{eq:quadbound}
\overline {\bf{Q}}_\phi = \left[ {\begin{array}{*{20}c}
   2 & 0 \\
   0 & 2 \\
\end{array}} \right].
\end{equation}

A scaling of $4 : 2 : 1 : 1.5$ is chosen for $g_{p,1} : g_{p,2} : g_{p,3} : \phi_p$, while $\overline {\boldsymbol \epsilon}$, $\boldsymbol{\overline \delta}_g$, and $\overline \delta_\phi$ are set to unity, and an approximation tolerance of $10^{-8}$ is used for $\underline {\boldsymbol \epsilon}$, $\boldsymbol{\underline \delta}_g$, and $\underline \delta_\phi$.

We note that the concave constraints serve a dual purpose. The first is the general role of uncertain hard constraints, while the second has more to do with their topological properties -- it is particularly difficult for an algorithm to remain feasible and avoid premature convergence while navigating around constraints with concave properties (see, e.g., the example in Figure \ref{fig:feasbreak}). Here, we construct a rather nasty case so as to highlight the benefits that enforcing the SCFO may bring to an algorithm -- this corresponds to the cases where the algorithms are initialized from ${\bf{u}}_0 = [-0.5, 0.05]$. To be more -- or, perhaps, less -- realistic, we also consider a nicer case with the initial point ${\bf{u}}_0 = [0, 0.4]$ where these constraints do not really play a role. We also note that their presence does not change the fact that this problem only has a single stable KKT point at ${\bf{u}}^* = [0.35, 0.32]$, which is, in this case, the global minimum. An unstable KKT point at ${\bf u} = [-0.09, 0.11]$ is also present, however.

As all three constraints are defined as hard constraints, we apply the filter gain criterion (\ref{eq:gainupper}) so as to ensure that the constraints are never violated. We also look at the cases where only Condition (\ref{eq:suff1}) is applied, as was done in Section 3. For these schemes, we do not auto-tune ${\boldsymbol \epsilon}$ and $\boldsymbol{\delta}_g$, simply setting them to ${\boldsymbol \epsilon} := 0.01 \overline {\boldsymbol \epsilon}$ and $\boldsymbol{\delta}_g := 0.01\boldsymbol{\overline \delta}_g$.

Five different RTO algorithms are considered, of which only two are reported here in the interest of space, with the results for the remaining three relegated to the Supplementary Material. We describe the two algorithms reported here below:

\vspace{2mm}
\noindent {\bf{Algorithm 1} -- Two-Step Approach}
\vspace{2mm}

A model of the plant, with a set of uncertain parameters $\boldsymbol{\theta}$, is available:

\begin{equation}\label{eq:ex4probPE}
\begin{array}{l}
\phi({\bf{u}},\boldsymbol{\theta}) = \theta_1 (u_1-0.3)^2 + \theta_2 (u_2-0.3)^2 \\
g_{1}({\bf{u}},\boldsymbol{\theta})= -\theta_3 u^2_1 - 3.5u_1 + u^2_2 -0.6 \le 0 \vspace{1mm} \\
g_{2}({\bf{u}},\boldsymbol{\theta})= \theta_4 u_1 + u_2 + \theta_5 \le 0 \vspace{1mm} \\
g_{3}({\bf{u}},\boldsymbol{\theta})= -u^2_1 - (u_2-0.15)^2 +0.01 \leq 0 \vspace{1mm} \\
\end{array},
\end{equation}

\noindent where we allow the third constraint to be modeled perfectly.

The parameters are re-estimated at each iteration (via linear regression) and the updated model is optimized numerically to compute ${\bf{u}}^*_{k+1}$. This represents a relatively good model-based RTO algorithm that nevertheless fails to achieve KKT convergence due to structural mismatch between the parametric model and the plant.

\vspace{2mm}
\noindent {\bf{Algorithm 2} -- Random Step}
\vspace{2mm}

Here, we draw the ${\bf{u}}^*_{k+1}$ out of a hat, using a uniform random number generator to provide an input such that $u^*_{k+1,1} \sim \mathcal{U}[-0.5, 0.5], u^*_{k+1,2} \sim \mathcal{U}[0, 0.8]$. While highly impractical and not at all advised, we use this method to underline the usefulness that enforcing the SCFO has, even when the RTO algorithm is completely haphazard.

\vspace{2mm}

\begin{figure}
\begin{center}
\includegraphics[width=12cm]{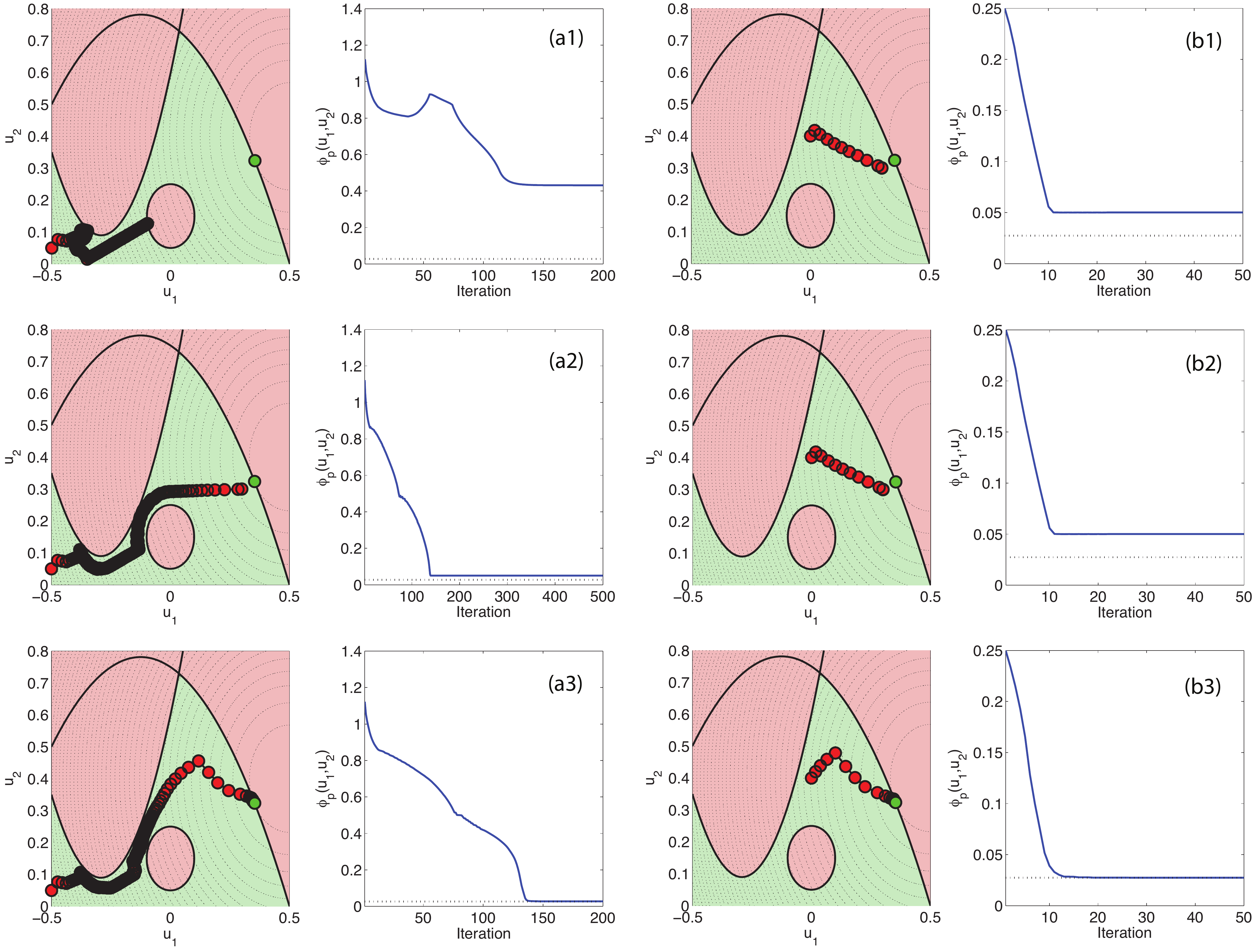}    % The printed column width  is 8.4 cm.
\caption{Performance of Algorithm 1 (two step) for the cases where (1) neither Conditions (\ref{eq:suff1}) nor (\ref{eq:suff2}) are enforced, (2) only Condition (\ref{eq:suff1}) is enforced, (3) the full SCFO are enforced, for (a) the first initial point and (b) the second initial point. Contours of the cost are given in black, with the optimal plant cost given by the dotted black lines in the right-hand figures. The plant optimum is plotted in green. We note that the SCFO are needed to enforce convergence to the optimum due to structural errors in the parametric model.}
\label{fig:ex4D}
\end{center}
\end{figure}

\begin{figure}
\begin{center}
\includegraphics[width=12cm]{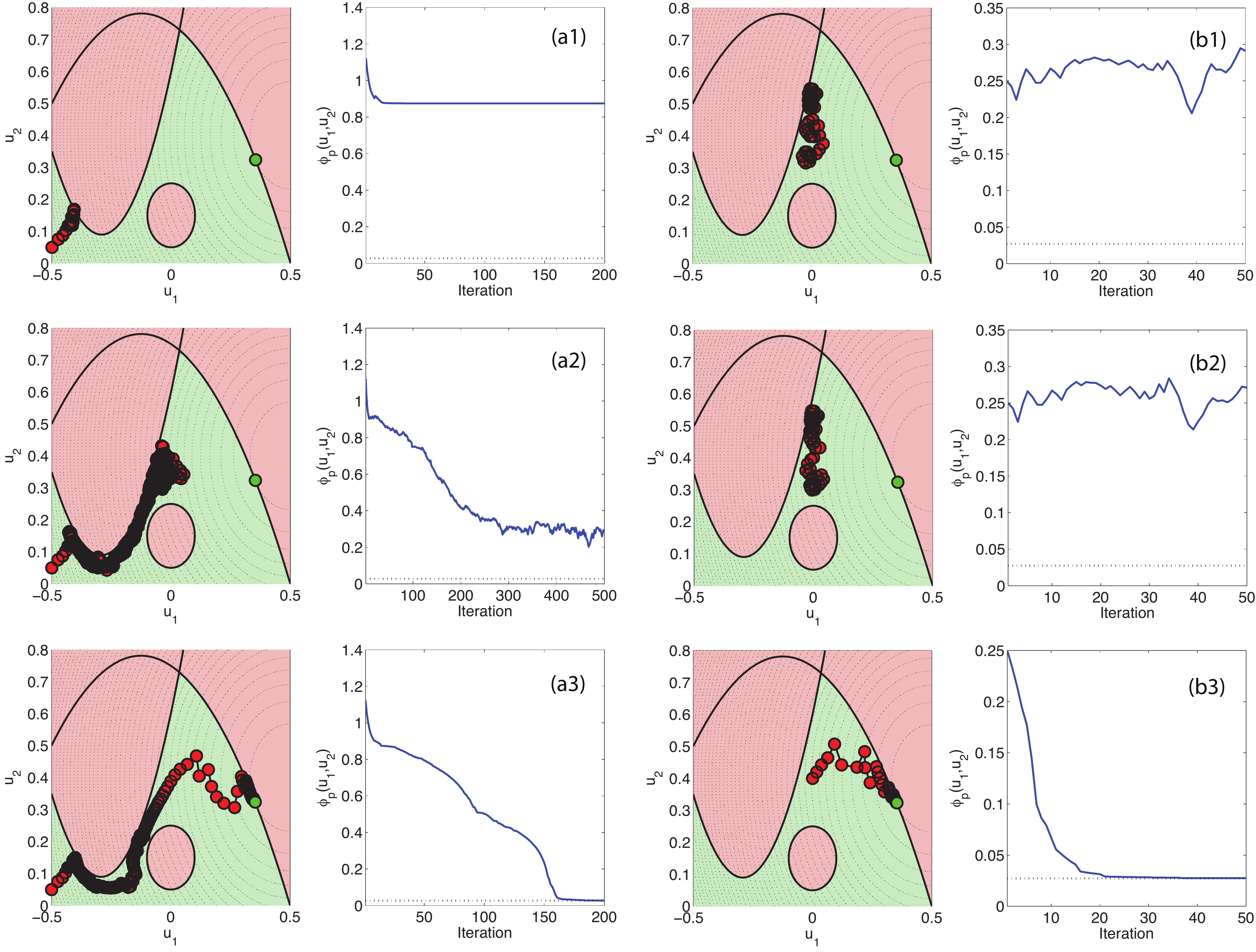}    % The printed column width  is 8.4 cm.
\caption{Performance of Algorithm 2 (random step) for the six cases. Not surprisingly, the algorithm does almost nothing desirable unless the SCFO are enforced.}
\label{fig:ex4E}
\end{center}
\end{figure}

Figures \ref{fig:ex4D}-\ref{fig:ex4E} provide the results. Leaving the algorithm-related remarks to the figure captions, we proceed to outline the noticeable general trends below:

\begin{itemize}
\item Feasibility is indeed preserved due to (\ref{eq:gainupper}) for all cases.
\item Premature convergence to one of the concave constraints is noted for all algorithms unless Condition (\ref{eq:suff1}) is enforced, in which case the iterates are able to ``slide'' around these obstacles.
\item Enforcing the SCFO leads to feasible convergence to the optimum in every case. This holds even when the algorithm is completely haphazard (Algorithm 2).
\item The SCFO always avoid the unstable KKT point at ${\bf u} = [-0.09, 0.11]$.
\item Enforcing Condition (\ref{eq:suff1}) alone proves sufficient for certain algorithms (See Algorithms A1, A2, and A3 in the Supplementary Material) to reach a region very close to the optimum. However, the lack of an auto-tuning scheme for the projection parameters keeps these realizations oscillating with an offset.
\item The steps are significantly smaller and progress significantly slower when the algorithm is close to a constraint. This is due to the numerator of Condition (\ref{eq:gainupper}).
\item Enforcing the SCFO for cases where the algorithm converges without them does affect performance slightly, in that more care is given with respect to the constraints, which affects the convergence trajectory (see examples in the Supplementary Material).
\item The choice of algorithm does not appear to be crucial when the SCFO are enforced, but does affect performance. For the case with the second initial point, for example, one may see that Algorithm 2 has reduced performance when compared to the others, due to its somewhat erratic trajectory.
\item Enforcing Condition (\ref{eq:suff2}) does indeed lead to monotonic improvement in the cost function.
\end{itemize}

\subsection{Summary}

The problem of convergence to a KKT optimum has been defined and the extents to which various RTO algorithm classes could guarantee this property have been reviewed, with the SCFO being proposed as a general set of conditions to force this property for any algorithm -- taking the feasibility conditions of the previous section and adding to them the requirement that the cost be locally decreasing at every iteration and that an additional constraint on the filter gain be satisfied to guarantee that the cost decrease for each step of the algorithm. Using Theorem 1, we showed that any algorithm meeting the SCFO would continue to decrease the cost until it became impossible to do so for a sufficiently low choice of projection parameters, which could only occur at a (practically stable) KKT point. This was demonstrated using several challenging examples.

\section{Conclusions}

\begin{figure}
\begin{center}
\includegraphics[width=6cm]{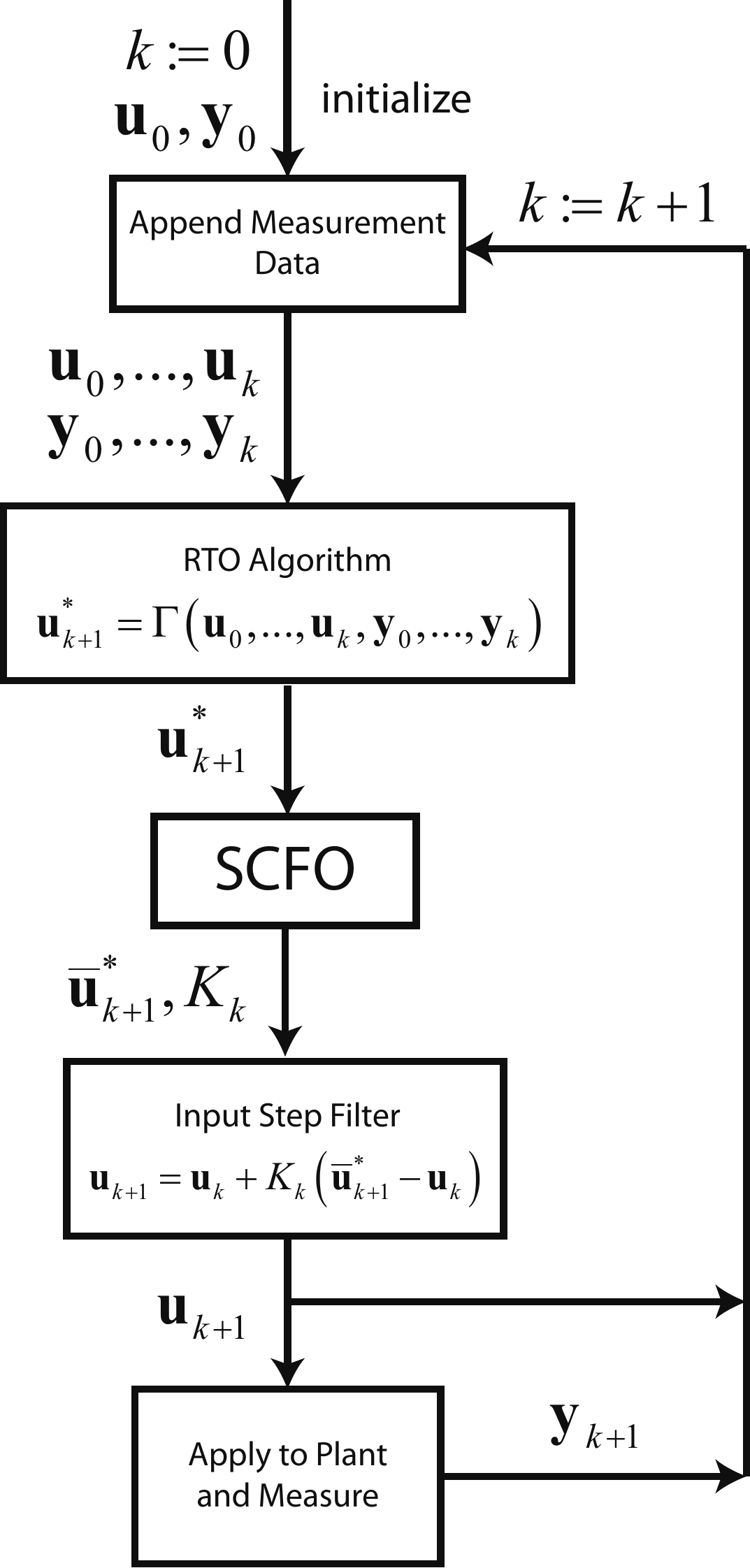}    % The printed column width  is 8.4 cm.
\caption{Generalized feedback structure of RTO schemes with the SCFO included so as to enforce feasibility and KKT convergence.}
\label{fig:fback2}
\end{center}
\end{figure}

The goal of the present paper has been to propose a set of sufficient conditions for convergence, via feasible iterates, to a plant KKT point in the context of the generalized RTO algorithm, and may be summarized via Figure \ref{fig:fback2} (as an extension of the structure in Figure \ref{fig:fback}). Sections 3 and 4 have presented the necessary theory with respect to feasibility and optimality, respectively, with the latter subsuming the former to comprise the full SCFO set. Testing the SCFO in simulation gave excellent performance, and showed that we could achieve convergence even for cases where the feasible set had unfavorable topological properties that made it impossible for most algorithms to get around unless the SCFO were applied.

It would be fair to note that, at least from a mathematical point of view, the proposed conditions are not terribly surprising, or even novel, in their structure -- simply put, they restate the somewhat natural fact that local derivative information, when coupled with higher-order global bounds, allow us to manipulate any optimization algorithm to force it to converge to a KKT point regardless of its innate structure. This sort of idea cannot possibly be new in the context of numerical optimization. However, we believe that the work in this paper constitutes an important contribution to the RTO field, where the numerous challenges, as outlined in the introduction, have led to the creation of many methods but no general set of guarantees. It is hoped that the SCFO can fill an important gap by providing a theoretical foundation to guide the development and operation of future (and current) RTO algorithms.

Another point worth noting is how the \emph{ad hoc} nature of the different algorithms is significantly reduced when placed into this framework. As we see from the proposed projection steps, the RTO algorithm still plays a role (by calculating the initial ${\bf{u}}^*_{k+1}$) in determining the general direction for the iterates to follow. As such, we naturally expect better, more accurate algorithms to get us to the optimum quicker, and are able to confirm these expectations in the simulated trials of Section 4. However, it is also important to note that the actual choice of algorithm is no longer crucial and is merely a preference -- when the SCFO are enforced, feasible-side convergence to a KKT point is ensured regardless. We believe that this has the potential to remove a considerable burden from the practitioner, as no ``optimal'' algorithm exists and several valid approaches could be proposed depending on the specific problem at hand.

Perhaps the immediate question introduced in this work is the following: should these conditions actually be enforced in practice as suggested in Figure \ref{fig:fback2}? The short and immediate answer is that they probably should, as with a relatively compact set of assumptions they allow for very strong and useful guarantees, as well as for a general foundation for RTO analysis. A longer answer might involve the study of the impact of the employed projection on the convergence rate of an algorithm, as it is more than likely that there exist problems for which such a projection actually makes convergence slower. In these latter cases, it may be better to let the algorithm operate as it would normally, and then enforce the SCFO only if it fails to perform as desired. We have not attempted to carry out such a study here.

The other question -- one of pressing importance and the subject of the companion paper -- is whether or not it is possible to actually enforce these conditions for a real system, as they depend on the knowledge of exact local derivatives, global upper bounds, and noise-free measurements, none of which is generally \emph{known} in application. It will be argued in the second paper that, while the SCFO cannot be applied as easily as they were here when the knowledge assumptions were made, they can still be robustly enforced for a number of practical realizations, and that doing so will generally lead to superior performance.

\bibliography{feasconv}             % bib file to produce the bibliography
                                                    % with bibtex (preferred)

\section*{Appendix}

\subsection*{Proof of Lemma 2}

We prove the lemma by considering the inequality along the line segment between an arbitrary pair ${\bf u}_k, {\bf u}_{k+1} \in \mathcal{I}$. The following one-dimensional parameterization is used:

\begin{equation}\label{eq:Qproof1}
f^u (\gamma) = f({\bf u} (\gamma)),
\end{equation}

\noindent with ${\bf u} (\gamma) = {\bf u}_k + \gamma ({\bf u}_{k+1} - {\bf u}_{k}), \; \gamma \in [0,1]$. As $f$ is twice continuously differentiable, it follows that $f^u$ is as well, which allows us to use the Taylor series expansion between $\gamma = 0$ and $\gamma = 1$, together with the mean-value theorem, to state \citep[Section 4.10-4]{Korn:00}:

\begin{equation}\label{eq:Qproof2}
\begin{array}{l}
f^u (1) = f^u (0) + \displaystyle \frac{df^u}{d \gamma} \Big |_{\gamma = 0} + R_1(1,0) \vspace{1mm} \\
R_1(1,0) = \displaystyle \frac{1}{2} \frac{d^2 f^u}{d\gamma^2} \Big |_{\gamma = \tilde \gamma}
\end{array},
\end{equation}

\noindent for some $\tilde \gamma \in (0,1)$. We proceed to define the first- and second-order derivatives in terms of the original function $f$. To do this we apply the chain rule:

\begin{equation}\label{eq:Qproof3}
\displaystyle \frac{df^u}{d \gamma} \Big |_{\gamma} = \displaystyle \mathop {\sum} \limits_{i = 1}^{n_u} \frac{\partial f}{\partial u_i} \Big |_{{\bf u} (\gamma)} \frac{d u_i}{d \gamma} \Big |_{\gamma} = \nabla f({\bf u}(\gamma))^T ({\bf u}_{k+1} - {\bf u}_{k}),
\end{equation}

\noindent and then differentiate once more with respect to $\gamma$:

\begin{equation}\label{eq:Qproof4}
\displaystyle \frac{d^2 f^u}{d \gamma^2} \Big |_{\gamma} = \displaystyle \mathop {\sum} \limits_{i = 1}^{n_u} \frac{d}{d \gamma} \left( \frac{\partial f}{\partial u_i} \Big |_{{\bf u} (\gamma)} \frac{d u_i}{d \gamma} \Big |_{\gamma} \right) = \displaystyle \mathop {\sum} \limits_{i = 1}^{n_u} \frac{d}{d \gamma} \left( \frac{\partial f}{\partial u_i} \Big |_{{\bf u} (\gamma)} \right) \frac{d u_i}{d \gamma} \Big |_{\gamma},
\end{equation}

\noindent where we have ignored the terms corresponding to $d^2 u_i / d \gamma^2$ as all such terms are 0. Applying the chain rule again yields:

\begin{equation}\label{eq:Qproof5}
\begin{array}{l}
\displaystyle \frac{d^2 f^u}{d \gamma^2} \Big |_{\gamma} = \displaystyle \mathop {\sum} \limits_{i = 1}^{n_u} \mathop {\sum} \limits_{j = 1}^{n_u}  \frac{\partial^2 f}{\partial u_i \partial u_j} \Big |_{{\bf u} (\gamma)} \frac{d u_j}{d \gamma} \Big |_{\gamma} \frac{d u_i}{d \gamma} \Big |_{\gamma} \\
= \displaystyle \mathop {\sum} \limits_{i = 1}^{n_u} \mathop {\sum} \limits_{j = 1}^{n_u}  \frac{\partial^2 f}{\partial u_i \partial u_j} \Big |_{{\bf u} (\gamma)} (u_{k+1,i} - u_{k,i})(u_{k+1,j} - u_{k,j})
\end{array}.
\end{equation}

Substituting the results of (\ref{eq:Qproof3}) and (\ref{eq:Qproof5}) into (\ref{eq:Qproof2}), and noting that $f^u (0) = f({\bf u}_k)$ and $f^u (1) = f({\bf u}_{k+1})$, leads to: 

\begin{equation}\label{eq:Qproof6}
\begin{array}{l}
f ({\bf u}_{k+1}) = f ({\bf u}_{k}) + \nabla f({\bf u}_{k})^T ({\bf u}_{k+1} - {\bf u}_{k}) + R_1(1,0) \vspace{2mm} \\
R_1(1,0) = \displaystyle \frac{1}{2} \mathop {\sum} \limits_{i = 1}^{n_u} \mathop {\sum} \limits_{j = 1}^{n_u}  \frac{\partial^2 f}{\partial u_i \partial u_j} \Big |_{{\bf u} (\tilde \gamma)} (u_{k+1,i} - u_{k,i})(u_{k+1,j} - u_{k,j})
\end{array}.
\end{equation}

We complete the proof by deriving an upper bound on the remainder term. First, note that (by $ab \leq |a||b|$):

\begin{equation}\label{eq:abs}
\begin{array}{l}
\displaystyle \frac{\partial^2 f}{\partial u_i \partial u_j} \Big |_{{\bf u} (\tilde \gamma)} (u_{k+1,i} - u_{k,i})(u_{k+1,j} - u_{k,j}) \vspace{1mm}\\
\leq \Bigg | \displaystyle \frac{\partial^2 f}{\partial u_i \partial u_j} \Big |_{{\bf u} (\tilde \gamma)} \Bigg | \Big | (u_{k+1,i} - u_{k,i})(u_{k+1,j} - u_{k,j}) \Big |
\end{array}.
\end{equation}

From the positivity of a quadratic (i.e. $0 \leq a^2 \pm 2ab + b^2 \Rightarrow 2|ab| \leq a^2 + b^2$), we have:

\begin{equation}\label{eq:Qproof7}
\begin{array}{l}
\Big | (u_{k+1,i} - u_{k,i})(u_{k+1,j} - u_{k,j}) \Big | \vspace{1mm}\\
\hspace{5mm}\leq \displaystyle \frac{1}{2} (u_{k+1,i} - u_{k,i})^2 + \frac{1}{2} (u_{k+1,j} - u_{k,j})^2
\end{array},
\end{equation}

\noindent from which it follows that:

\begin{equation}\label{eq:Qproof8}
\begin{array}{l}
\displaystyle \frac{\partial^2 f}{\partial u_i \partial u_j} \Big |_{{\bf u} (\tilde \gamma)} (u_{k+1,i} - u_{k,i})(u_{k+1,j} - u_{k,j}) \vspace{1mm}\\
\leq \displaystyle \frac{1}{2} \Bigg | \frac{\partial^2 f}{\partial u_i \partial u_j} \Big |_{{\bf u} (\tilde \gamma)} \Bigg | \left[ (u_{k+1,i} - u_{k,i})^2 + (u_{k+1,j} - u_{k,j})^2 \right] 
\end{array}.
\end{equation}

Substituting this into (\ref{eq:Qproof6}) then yields the following bound on the remainder:

\begin{equation}\label{eq:Qproof9}
\begin{array}{l}
R_1(1,0) \leq \displaystyle \frac{1}{4}  \displaystyle \mathop {\sum} \limits_{i = 1}^{n_u} \mathop {\sum} \limits_{j = 1}^{n_u}  \Bigg | \frac{\partial^2 f}{\partial u_i \partial u_j} \Big |_{{\bf u} (\tilde \gamma)} \Bigg |  (u_{k+1,i} - u_{k,i})^2 \\
\hspace{14mm}+ \displaystyle \frac{1}{4}  \displaystyle \mathop {\sum} \limits_{i = 1}^{n_u} \mathop {\sum} \limits_{j = 1}^{n_u}  \Bigg | \frac{\partial^2 f}{\partial u_i \partial u_j} \Big |_{{\bf u} (\tilde \gamma)} \Bigg |  (u_{k+1,j} - u_{k,j})^2
\end{array}.
\end{equation}

By Clairaut's theorem, we have that:

\begin{equation}\label{eq:clairaut}
\frac{\partial^2 f}{\partial u_i \partial u_j} \Big |_{{\bf u} (\tilde \gamma)} = \frac{\partial^2 f}{\partial u_j \partial u_i} \Big |_{{\bf u} (\tilde \gamma)},
\end{equation}

\noindent which, together with the interchangeability of the order of summation \citep[Section 4.8-3]{Korn:00}, allows us to rewrite the second term on the right-hand side of (\ref{eq:Qproof9}) as:

\begin{equation}\label{eq:inter}
\frac{1}{4}  \displaystyle \mathop {\sum} \limits_{j = 1}^{n_u} \mathop {\sum} \limits_{i = 1}^{n_u}  \Bigg | \frac{\partial^2 f}{\partial u_j \partial u_i} \Big |_{{\bf u} (\tilde \gamma)} \Bigg |  (u_{k+1,j} - u_{k,j})^2,
\end{equation}

\noindent which is clearly equivalent to the first term (only the choice of indices differs). This allows us to combine the two to obtain:

\begin{equation}\label{eq:Qproof10}
R_1(1,0) \leq \frac{1}{2}  \displaystyle \mathop {\sum} \limits_{i = 1}^{n_u} \mathop {\sum} \limits_{j = 1}^{n_u}  \Bigg | \frac{\partial^2 f}{\partial u_i \partial u_j} \Big |_{{\bf u} (\tilde \gamma)} \Bigg |  (u_{k+1,i} - u_{k,i})^2.
\end{equation}

To make this bound global for any choice of ${\bf u}_k, {\bf u}_{k+1} \in \mathcal{I}$ and all $\tilde \gamma \in (0,1)$, we use (\ref{eq:M0}), which implies: 

\begin{equation}\label{eq:Qproof11}
\Bigg | \frac{\partial^2 f}{\partial u_i \partial u_j} \Big |_{{\bf u}} \Bigg |  < M_{ij}, \;\; \forall {\bf u} \in {\mathcal{I}}.
\end{equation}

\noindent This globally bounds the remainder as:

\begin{equation}\label{eq:Qproof12}
R_1(1,0) \leq \frac{1}{2}  \displaystyle \mathop {\sum} \limits_{i = 1}^{n_u} \mathop {\sum} \limits_{j = 1}^{n_u}  M_{ij}  (u_{k+1,i} - u_{k,i})^2,
\end{equation}

\noindent or, in vector notation, as:

\begin{equation}\label{eq:Qproof13}
R_1(1,0) \leq \frac{1}{2} ({\bf u}_{k+1} - {\bf u}_k)^T \overline {\bf Q} ({\bf u}_{k+1} - {\bf u}_k),
\end{equation}

\noindent with $\overline {\bf Q}$ a diagonal matrix with the diagonal elements defined as in (\ref{eq:Qbounddiag}). Substituting (\ref{eq:Qproof13}) into (\ref{eq:Qproof6}) and rearranging leads to the main result. \hfill $\square$

\subsection*{Proof of Theorem 3}

We first obtain a more global version of (\ref{eq:gainupper}) that is independent of the RTO algorithm by remarking that:

\begin{equation}\label{eq:cor1p3}
\displaystyle \sum\limits_{i = 1}^{n_u} {\kappa_{ji} | u^*_{k+1,i} - u_{k,i} | } \leq {\boldsymbol \kappa}_j^T \left( {\bf u}^U - {\bf u}^L \right), \; \forall j = 1,...,n_g,
\end{equation} 
 
\noindent which allows us to lower bound $K_k$ as:

\begin{equation}\label{eq:Klimit0}
0 < \mathop {\min} \limits_{j = 1,...,n_g} \left[ \frac{ -g_{p,j}({\bf{u}}_{k } )  }{{\boldsymbol \kappa}_j^T ({\bf u}^U - {\bf u}^L)} \right] \leq  \mathop{\min}\limits_{j = 1,...,n_g} \left[ \frac{{-g_{p,j}({\bf{u}}_{k } )}}{\displaystyle \sum\limits_{i = 1}^{n_u} {\kappa_{ji} | u^*_{k+1,i} - u_{k,i} | }} \right] = K_k.
\end{equation} 

\noindent This then leaves us with the simpler task of lower bounding $K_k$ by lower bounding $-g_{p,j}({\bf u}_k)$, since ${\boldsymbol \kappa}_j^T ({\bf u}^U - {\bf u}^L)$ is independent of iteration.

We proceed by stating that $-g_{p,j}({\bf u}_k)$ may be lower bounded differently depending on its value. Specifically, we break the full set of possibilities into the following subcases:

\begin{enumerate}
\item $-g_{p,j}({\bf u}_k) > \epsilon_j$
\item $-g_{p,j}({\bf u}_k) \leq \epsilon_j$
\begin{enumerate}
\item $-g_{p,j}({\bf{u}}_{k } ) \geq \displaystyle \frac{K_{\epsilon,j} \delta_{g,j}}{\gamma_\kappa}$
\item $-g_{p,j}({\bf{u}}_{k } )  < \displaystyle \frac{K_{\epsilon,j} \delta_{g,j}}{\gamma_\kappa}$
\end{enumerate}
\end{enumerate}

\noindent where we will refer to the three subcases as Case 1, Case 2a, and Case 2b. We will now derive the lowest values that $-g_{p,j}$ can achieve for each of these three scenarios.

For Case 1, we start by noting that, from the definition of $\gamma_\kappa$ in (\ref{eq:gammakappa0}), we have that:

\begin{equation}\label{eq:gammakappa}
\tilde \kappa_{ji} \leq \gamma_\kappa \kappa_{ji}, \; \forall i = 1,...,n_u, \;\; j = 1,...,n_g,
\end{equation}

\noindent with $\gamma_\kappa \in (0,1)$ implicit from the strictness of $\kappa$ and the fact that all of the $\tilde \kappa$ cannot be null simultaneously (each constraint being a function of the inputs and not simply a constant value). We now follow the same steps as in Lemma 1 and Theorem 2, this time for the case with the nonstrict Lipschitz constants, to obtain:

\begin{equation}\label{eq:cor1p1}
\begin{array}{l}
g_{p,j} ({\bf{u}}_{k+1}) \leq g_{p,j} ({\bf{u}}_{k}) + \displaystyle \sum\limits_{i = 1}^{n_u} {\tilde \kappa_{ji} | u_{k+1,i} - u_{k,i} | } \\
\Rightarrow g_{p,j} ({\bf{u}}_{k+1}) \leq g_{p,j} ({\bf{u}}_{k}) + \gamma_\kappa \displaystyle \sum\limits_{i = 1}^{n_u} { \kappa_{ji} | u_{k+1,i} - u_{k,i} | }
\end{array}.
\end{equation}

Since

\begin{equation}\label{eq:cor1p2}
\begin{array}{l}
g_{p,j} ({\bf{u}}_{k}) + \displaystyle \sum\limits_{i = 1}^{n_u} {\kappa_{ji} | u_{k+1,i} - u_{k,i} | } \leq 0 \vspace{1mm} \\
\hspace{10mm}\Leftrightarrow \displaystyle \sum\limits_{i = 1}^{n_u} {\kappa_{ji} | u_{k+1,i} - u_{k,i} | } \leq -g_{p,j} ({\bf u}_k)
\end{array}
\end{equation}

\noindent is enforced by (\ref{eq:gainupper}) at all iterations, it follows that:

\begin{equation}\label{eq:cor1p3}
\begin{array}{l}
g_{p,j} ({\bf{u}}_{k+1}) \leq g_{p,j} ({\bf{u}}_{k}) - \gamma_\kappa g_{p,j} ({\bf u}_k) = (1- \gamma_\kappa) g_{p,j} ({\bf u}_k) \vspace{1mm} \\
\hspace{10mm}\Leftrightarrow -g_{p,j} ({\bf{u}}_{k+1}) \geq (1- \gamma_\kappa) (-g_{p,j} ({\bf u}_k))
\end{array}.
\end{equation}

\noindent It is evident that the lowest value that $-g_{p,j} ({\bf{u}}_{k+1})$ can achieve while $-g_{p,j} ({\bf{u}}_{k}) > \epsilon_j$ is:

\begin{equation}\label{eq:cor1p8}
-g_{p,j}({\bf{u}}_{k +1 } ) > (1-\gamma_\kappa) \epsilon_j,
\end{equation}

\noindent after which Case 1 would no longer be pertinent and we would shift our analysis to Cases 2a and 2b.

For these cases, we employ the result of Lemma 3, which states that:

\begin{equation}\label{eq:descentstep1}
0 < K_k < -2\frac{\nabla g_{p,j}({\bf{u}}_k)^T ({\bf{u}}^*_{k+1} - {\bf{u}}_{k})}{({\bf{u}}^*_{k+1} - {\bf{u}}_{k})^T \overline {\bf{Q}}_j ({\bf{u}}^*_{k+1} - {\bf{u}}_{k})} \Rightarrow -g_{p,j}({\bf{u}}_{k+1}) > -g_{p,j}({\bf{u}}_k).
\end{equation}

Using Condition (\ref{eq:suff1}), together with the fact that $\overline {\bf Q}_j \succ 0$, allows:

\begin{equation}\label{eq:Keps0}
K_{\epsilon,j} \leq -2\frac{\nabla g_{p,j}({\bf{u}}_k)^T ({\bf{u}}^*_{k+1} - {\bf{u}}_{k})}{({\bf{u}}^*_{k+1} - {\bf{u}}_{k})^T \overline {\bf{Q}}_j ({\bf{u}}^*_{k+1} - {\bf{u}}_{k})},
\end{equation}

\noindent as (\ref{eq:Keps}) represents the minimal value that the right-hand side could possibly take for any iteration where $g_{p,j}$ is $\epsilon$-active -- this is done by maximizing the (negative) numerator and maximizing the (positive) denominator. This allows us to make (\ref{eq:descentstep1}) independent of the RTO algorithm and to state:

\begin{equation}\label{eq:Kepsbar}
K_k < K_{\epsilon,j} \Rightarrow -g_{p,j}({\bf{u}}_{k+1}) > -g_{p,j}({\bf{u}}_k).
\end{equation} 

\noindent Using (\ref{eq:Klimit0}), we may extend this statement further:

\begin{equation}\label{eq:Kepsbar2}
\mathop {\min} \limits_{\hat j = 1,...,n_g} \left[ \frac{ -g_{p,\hat j}({\bf{u}}_{k } )  }{\displaystyle \sum\limits_{i = 1}^{n_u} {\kappa_{\hat j i} | u^*_{k+1,i} - u_{k,i} | }} \right] < K_{\epsilon,j} \Rightarrow -g_{p,j}({\bf{u}}_{k+1}) > -g_{p,j}({\bf{u}}_k),
\end{equation} 

\noindent where we make a temporary change of indices so as to distinguish between the terms in the minimum operator (indexed by $\hat j$) from the single constraint whose evolution is being considered (indexed by $j$). However, since satisfying the inequality for any arbitrary $j$ also satisfies it for the minimum:

\begin{equation}\label{eq:imp1}
\frac{ -g_{p,j}({\bf{u}}_{k } )  }{\displaystyle \sum\limits_{i = 1}^{n_u} {\kappa_{ji} | u^*_{k+1,i} - u_{k,i} | }} < K_{\epsilon,j} \Rightarrow \mathop {\min} \limits_{\hat j = 1,...,n_g} \left[ \frac{ -g_{p,\hat j}({\bf{u}}_{k } )  }{\displaystyle \sum\limits_{i = 1}^{n_u} {\kappa_{\hat j i} | u^*_{k+1,i} - u_{k,i} | }} \right] < K_{\epsilon,j},
\end{equation} 

\noindent we may simplify (\ref{eq:Kepsbar2}) to:

\begin{equation}\label{eq:Kepsbar3}
\frac{ -g_{p,j}({\bf{u}}_{k } )  }{\displaystyle \sum\limits_{i = 1}^{n_u} {\kappa_{j i} | u^*_{k+1,i} - u_{k,i} | }} < K_{\epsilon,j} \Rightarrow -g_{p,j}({\bf{u}}_{k+1}) > -g_{p,j}({\bf{u}}_k).
\end{equation} 

In order to make this statement independent of the RTO algorithm, we exploit the $\epsilon$-activity of $g_{p,j}$ and rewrite the condition in (\ref{eq:suff1}) as:

\begin{equation}\label{eq:suff1re}
\nabla g_{p,j}({\bf{u}}_k)^T ({\bf{u}}_{k} - {\bf{u}}^*_{k+1}) = \displaystyle \sum_{i=1}^{n_u} \frac{\partial g_{p,j}}{\partial u_i} \Big |_{{\bf u}_k} \left( u_{k,i} - u_{k+1,i}^* \right) \geq \delta_{g,j},
\end{equation}

\noindent from which it readily follows that:

\begin{equation}\label{eq:suff1re2}
\begin{array}{l}
\displaystyle \sum_{i=1}^{n_u} \tilde \kappa_{ji} | u_{k+1,i}^* - u_{k,i} | \geq \displaystyle \sum_{i=1}^{n_u} \frac{\partial g_{p,j}}{\partial u_i} \Big |_{{\bf u}_k} \left( u_{k,i} - u_{k+1,i}^* \right) \geq \delta_{g,j} \\
\Rightarrow \gamma_\kappa \displaystyle \sum_{i=1}^{n_u} \kappa_{ji} | u_{k+1,i}^* - u_{k,i} | \geq \delta_{g,j} \Leftrightarrow \displaystyle \sum_{i=1}^{n_u} \kappa_{ji} | u_{k+1,i}^* - u_{k,i} | \geq \frac{\delta_{g,j}}{\gamma_\kappa}
\end{array},
\end{equation}

\noindent thereby allowing us to further chain the implication of (\ref{eq:imp1}):

\begin{equation}\label{eq:imp2}
\frac{ -\gamma_\kappa g_{p,j}({\bf{u}}_{k } )  }{ {\delta_{g,j} }} < K_{\epsilon,j} \Rightarrow \frac{ -g_{p,j}({\bf{u}}_{k } )  }{\displaystyle \sum\limits_{i = 1}^{n_u} {\kappa_{ji} | u^*_{k+1,i} - u_{k,i} | }} < K_{\epsilon,j},
\end{equation} 

\noindent to obtain:

\begin{equation}\label{eq:Kepsbar4}
\begin{array}{l}
\displaystyle \frac{ -\gamma_\kappa g_{p,j}({\bf{u}}_{k } )  }{ {\delta_{g,j} }} < K_{\epsilon,j} \Rightarrow -g_{p,j}({\bf{u}}_{k+1}) > -g_{p,j}({\bf{u}}_k) \\
\Leftrightarrow -g_{p,j}({\bf{u}}_{k } )  < \displaystyle \frac{K_{\epsilon,j} \delta_{g,j}}{\gamma_\kappa} \Rightarrow -g_{p,j}({\bf{u}}_{k+1}) > -g_{p,j}({\bf{u}}_k)
\end{array}.
\end{equation} 

We are now equipped to make statements about Cases 2a and 2b. Note first that Case 2a may not exist -- this occurs if 

\begin{equation}\label{eq:casenull}
\epsilon_j  < \displaystyle \frac{K_{\epsilon,j} \delta_{g,j}}{\gamma_\kappa},
\end{equation} 

\noindent since Case 1 would remain dominant. If this is so, then only Case 1 and Case 2b need consideration. Assuming that Case 2a can occur, it follows that we can exploit the bound in (\ref{eq:cor1p3}) to calculate the lowest value that $-g_{p,j} ({\bf{u}}_{k+1})$ can achieve while $-g_{p,j} ({\bf{u}}_{k}) > \displaystyle \frac{K_{\epsilon,j} \delta_{g,j}}{\gamma_\kappa}$:

\begin{equation}\label{eq:cor1p7}
-g_{p,j}({\bf{u}}_{k+1 } ) > (1-\gamma_\kappa) \frac{K_{\epsilon,j} \delta_{g,j}}{\gamma_\kappa},
\end{equation}

\noindent after which we would shift to Case 2b.

In this final scenario, we have, by (\ref{eq:Kepsbar4}), that the value cannot decrease and that $-g_{p,j}({\bf{u}}_{k+1}) > -g_{p,j}({\bf{u}}_k)$. It thus only remains to lower bound $-g_{p,j}({\bf{u}}_k)$ in this bound to obtain the ultimate lower bound on $-g_{p,j}$ for all iterations. To do this, one needs to analyze how one can arrive at Case 2b. The first means is directly from Case 1, in which case we may shift the indices (by applying $k := k+1$) and note that:

\begin{equation}\label{eq:cor1p8a}
-g_{p,j}({\bf{u}}_{k} ) > (1-\gamma_\kappa) \epsilon_j.
\end{equation}

The second means is to go to Case 2b from Case 2a, where, using the same logic, we have:

\begin{equation}\label{eq:cor1p7a}
-g_{p,j}({\bf{u}}_k ) > (1-\gamma_\kappa) \frac{K_{\epsilon,j} \delta_{g,j}}{\gamma_\kappa}.
\end{equation}

Finally, there is also the possibility that the RTO algorithm is initialized directly in Case 2b, with $-g_{p,j}({\bf{u}}_0 ) < \displaystyle  \frac{K_{\epsilon,j} \delta_{g,j}}{\gamma_\kappa}$. In this case, we have:

\begin{equation}\label{eq:cor1p7a}
-g_{p,j}({\bf{u}}_k ) > -g_{p,j}({\bf{u}}_0 ).
\end{equation}

To account for all of these possibilities, which are comprehensive and account for all that could occur, we may proceed to obtain the global lower bound on $-g_{p,j}$ by taking the minimum of the three:

\begin{equation}\label{eq:cor1p9}
-g_{p,j}({\bf{u}}_{k } ) >  \mathop {\min} \left( (1-\gamma_\kappa) \epsilon_j, (1-\gamma_\kappa) \frac{K_{\epsilon,j} \delta_{g,j}}{\gamma_\kappa}, -g_{p,j}({\bf{u}}_0 ) \right).
\end{equation}

This allows us to restate (\ref{eq:Klimit0}) as:

\begin{equation}\label{eq:Klimit1}
\begin{array}{l}
\displaystyle 0 < \mathop {\min} \limits_{j = 1,...,n_g} \left[ \frac{ \mathop {\min} \left( (1-\gamma_\kappa) \epsilon_j, (1-\gamma_\kappa) \displaystyle \frac{K_{\epsilon,j} \delta_{g,j}}{\gamma_\kappa}, -g_{p,j}({\bf{u}}_0 ) \right)  }{{\boldsymbol \kappa}_j^T ({\bf u}^U - {\bf u}^L)} \right] \vspace{1mm} \\
\hspace{45mm}\displaystyle < \mathop {\min} \limits_{j = 1,...,n_g} \left[ \frac{ -g_{p,j}({\bf{u}}_{k } )  }{{\boldsymbol \kappa}_j^T ({\bf u}^U - {\bf u}^L)} \right] \leq  K_k
\end{array},
\end{equation} 

\noindent which proves the desired result. \hfill $\square$

\subsection*{Proof of Theorem 5}

\begin{enumerate}[(i)]
\item (\ref{eq:gainupper}) and (\ref{eq:suff1}) guarantee the existence of $K_{min} > 0$, by Theorem 3, which in turn allows the application of (\ref{eq:suff2}) and the guarantee that the conditions may only be satisfied for a finite number of iterations, as upper bounded in (\ref{eq:upperiter}). As the algorithm is entirely deterministic, it follows that the inputs remain at ${\bf u}_\infty$ indefinitely following the first instance where the conditions can no longer be satisfied.
\item Feasibility follows from (\ref{eq:gainupper}).
\item Monotonic cost improvement until ${\bf u}_\infty$ follows from (\ref{eq:suff2}).
\item We must characterize the point ${\bf u}_\infty$ where the three sets of Conditions (\ref{eq:gainupper}), (\ref{eq:suff1}), and (\ref{eq:suff2}) cannot be satisfied. First, we drop the conditions on the filter gain (Condition (\ref{eq:gainupper}) and the second condition in (\ref{eq:suff2})) from the discussion, noting that these can always be satisfied for $K_k = K_{min}$ and are not a concern. We thus turn to the inequality constraints of Conditions (\ref{eq:suff1}) and (\ref{eq:suff2}), as well as to the box constraints. Using the same notation as in Theorem 4, these may be rewritten as:

\begin{equation}\label{eq:optre}
\begin{array}{l}
\left[ {\begin{array}{*{20}c}
   \nabla \phi_p ({\bf{u}}_\infty)^T  \\
   {\bf{J}}_\infty^T \\
\end{array}} \right]({\bf{u}} - {\bf{u}}_{\infty}) \preceq \left[ {\begin{array}{*{20}c}
   -\delta_\phi \\
   -\boldsymbol{\delta}_g  \\
   {\bf{0}} \\
\end{array}} \right] \vspace{1mm} \\
\tilde {\bf{u}}^L \preceq {\bf{u}} \preceq \tilde {\bf{u}}^U
\end{array}.
\end{equation}

Equivalently, we may consider the scaled version:

\begin{equation}\label{eq:optrescaled}
\begin{array}{l}
\left[ {\begin{array}{*{20}c}
   \nabla \phi_p ({\bf{u}}_\infty)^T  \\
   {\bf{J}}_\infty^T \\
\end{array}} \right] \alpha ({\bf{u}} - {\bf{u}}_{\infty}) \preceq \left[ {\begin{array}{*{20}c}
   -\alpha \delta_\phi \\
   -\alpha \boldsymbol{\delta}_g  \\
   {\bf{0}} \\
\end{array}} \right] \vspace{1mm} \\
\tilde {\bf{u}}^L \preceq {\bf{u}}_\infty + \alpha ({\bf{u}} - {\bf{u}}_\infty) \preceq \tilde {\bf{u}}^U
\end{array}.
\end{equation}

Since $\tilde {\bf{u}}^L \prec {\bf{u}}_\infty \prec \tilde {\bf{u}}^U$, it follows that there exists a choice $\alpha^L > 0$, for any ${\bf u}$, so that the inactive box constraints will be satisfied. Since $\alpha \rightarrow 0 \Rightarrow \alpha \boldsymbol{\delta}_g, \alpha \delta_\phi \rightarrow {\bf 0}$, it follows that there exist small enough $\boldsymbol{\delta}^L_g = \alpha^L \boldsymbol{\delta}_g \succ {\bf 0}$ and $\delta_\phi^L = \alpha^L \delta_\phi > 0$ so that the inactive box constraints may be satisfied in the limit for this choice of $\boldsymbol{\delta}_g$ and $\delta_\phi$. We have therefore proven that the inactive box constraints cannot cause infeasibility in the conditions as $\boldsymbol{\delta}_g, \delta_\phi \rightarrow {\bf 0}$, and proceed to consider the infeasibility of the simpler case:

\begin{equation}\label{eq:optresimp}
\left[ {\begin{array}{*{20}c}
   \nabla \phi_p ({\bf{u}}_\infty)^T  \\
   {\bf{J}}_\infty^T \\
\end{array}} \right] \alpha^L ({\bf{u}} - {\bf{u}}_{\infty}) \preceq \left[ {\begin{array}{*{20}c}
   -\delta_\phi^L \\
   -\boldsymbol{\delta}_g^L  \\
   {\bf{0}} \\
\end{array}} \right],
\end{equation}

which is equivalent with respect to feasibility to the original unscaled set:

\begin{equation}\label{eq:optresimp2}
\left[ {\begin{array}{*{20}c}
   \nabla \phi_p ({\bf{u}}_\infty)^T  \\
   {\bf{J}}_\infty^T \\
\end{array}} \right] ({\bf{u}} - {\bf{u}}_{\infty}) \preceq \left[ {\begin{array}{*{20}c}
   -\delta_\phi \\
   -\boldsymbol{\delta}_g  \\
   {\bf{0}} \\
\end{array}} \right].
\end{equation}

We analyze the infeasibility of this set as ${\boldsymbol \epsilon} \rightarrow {\bf 0}$. From Theorem 1, we know that the system in (\ref{eq:optresimp2}) is infeasible iff:

\begin{equation}\label{eq:KKTstat}
a_\phi \nabla \phi_p ({\bf{u}}_\infty) + a_1 {\bf{J}}_{\infty,1} + ... + a_{n_J}{\bf{J}}_{\infty,n_J} = {\bf{0}},
\end{equation}

\noindent with at least one $a$ coefficient corresponding to either the cost or the uncertain constraints strictly positive.

Since (\ref{eq:KKTstat}) must be true at ${\bf u}_\infty$ for ${\bf 0} \prec \boldsymbol{\delta}_g \preceq \boldsymbol{\delta}_g^L, 0 < \delta_\phi \leq \delta_\phi^L$, we proceed to analyze the several cases that satisfy (\ref{eq:KKTstat}). The first corresponds to the occurrence of a trivial negative spanning where one of the elements in the summation is ${\bf 0}$. This may occur with either the cost or the uncertain constraints, as the box constraints cannot have a zero gradient. If this is true for the cost, with $\nabla \phi_p ({\bf{u}}_\infty) = {\bf 0}$, then the KKT conditions for an unconstrained KKT point are satisfied and the KKT error is 0. For the constraints, we have, from Assumption A3, that there exists $\epsilon_{m,j} > 0$ for every uncertain constraint $g_{p,j}$ so that the gradient $\nabla g_{p,j} ({\bf u}_\infty)$ cannot be ${\bf 0}$ as ${\boldsymbol \epsilon} \rightarrow {\bf 0}$.

We consider the cases of nontrivial negative spanning next, where at least two of the $a$ coefficients are strictly positive. There are two cases of interest -- one where $a_\phi = 0$ (the cost descent condition does not contribute to infeasibility) and one where $a_\phi > 0$ (the cost does contribute). Suppose first that $a_\phi = 0$ as ${\boldsymbol \epsilon} \rightarrow {\bf 0}$. If this is true, then the relevant constraint gradients span each other negatively even as the $\epsilon$-active set becomes an arbitrarily good approximation of the true active set. This, in turn, means that there is no direction that locally decreases all of these constraints, i.e. a locally feasible direction. This then implies that ${\bf u}_\infty$ is a singleton in the set defined by the active constraints, thereby allowing us to ignore this case. As stated before, such problems are ill-posed and do not require RTO. We are therefore left with the case of $a_\phi > 0$ for a sufficiently small ${\boldsymbol \epsilon}^L \succ {\bf 0}$ (we need ${\boldsymbol \epsilon}^L$ to be small enough so as to ensure that there is no negative spanning between the $\epsilon$-active constraints).

We now proceed to derive an upper bound on the KKT error that is valid for ${\bf 0} \prec {\boldsymbol \epsilon} \preceq {\boldsymbol \epsilon}^L$, ${\bf 0} \prec \boldsymbol{\delta}_g \preceq \boldsymbol{\delta}_g^L$, and $0 < \delta_\phi \leq \delta_\phi^L$ and show how it must tend to 0 as ${\boldsymbol \epsilon} \rightarrow {\bf 0}$.

For this case, we may, without loss of generality, scale to set $a_\phi = 1$ and, changing the notation on the coefficients, write (\ref{eq:KKTstat}) in a form that is analogous to that of the Lagrangian in (\ref{eq:KKT}):

\begin{equation}\label{eq:lagana}
\nabla \tilde{\mathcal{L}}({\bf u}_\infty) = \nabla \phi_p ({\bf u}_\infty) + \displaystyle \sum\limits_{j = 1}^{n_g}  \tilde \mu_j \nabla g_{p,j} ({\bf u}_\infty) -  \tilde{\boldsymbol{\zeta}}^L +  \tilde{\boldsymbol{\zeta}}^U  = {\bf 0},
\end{equation}

with $\tilde \mu_j = 0, \; \forall j : g_{p,j} ({\bf u}_\infty) < -\epsilon_j$, $\tilde \zeta_i^U = 0, \; \forall i : u_{\infty,i} < u_{i}^U$, $\tilde \zeta_i^L = 0, \; \forall i : u_{\infty,i} > u_{i}^L$ establishing equivalence with (\ref{eq:KKTstat}).

Using (\ref{eq:lagana}), we may express the gradient of the true Lagrangian as:

\begin{equation}\label{eq:lagref}
\begin{array}{l}
\nabla \mathcal{L}({\bf u}_\infty) = \nabla \mathcal{L}({\bf u}_\infty) - \nabla \tilde{\mathcal{L}}({\bf u}_\infty) =  \vspace{1mm} \\
\hspace{20mm}\displaystyle \sum\limits_{j = 1}^{n_g} (\mu_j - \tilde \mu_j) \nabla g_{p,j} ({\bf u}_\infty) -  (\boldsymbol{\zeta}^L - \tilde{\boldsymbol{\zeta}}^L) +  (\boldsymbol{\zeta}^U - \tilde{\boldsymbol{\zeta}}^U)
\end{array}.
\end{equation}

\noindent Noting that $\nabla \mathcal{L}({\bf u}_\infty) = {\bf 0}$ for the choice $ {\boldsymbol \mu} = \tilde {\boldsymbol \mu}$, $\boldsymbol{\zeta}^L = \tilde{\boldsymbol{\zeta}}^L$, $\boldsymbol{\zeta}^U = \tilde{\boldsymbol{\zeta}}^U$, we upper bound the KKT error as follows:

\begin{equation}\label{eq:KKTerrup}
\begin{array}{l}
\mathop {\inf} \limits_{{\boldsymbol \mu}, {\boldsymbol \zeta^L}, {\boldsymbol \zeta^U} \succeq {\bf 0}} \Big ( \nabla \mathcal{L}({\bf u}_\infty)^T \nabla \mathcal{L}({\bf u}_\infty) +  \displaystyle \sum\limits_{j = 1}^{n_g} \left[ \mu_j g_{p,j} ({\bf u}_\infty) \right]^2 + \\
\hspace{20mm}\displaystyle \sum\limits_{i = 1}^{n_u} \left[ ( \zeta^L_i (u^L_i - u_{\infty,i}) )^2 + ( \zeta^U_i (u_{\infty,i} - u^U_i) )^2  \right] \Big ) \\ 
\leq \mathop {\inf} \limits_{{\boldsymbol \mu} = \tilde{\boldsymbol \mu}, {\boldsymbol \zeta^L} = \tilde{\boldsymbol \zeta}^L, {\boldsymbol \zeta^U} = \tilde{\boldsymbol \zeta}^U} \Big ( \nabla \mathcal{L}({\bf u}_\infty)^T \nabla \mathcal{L}({\bf u}_\infty) +  \displaystyle \sum\limits_{j = 1}^{n_g} \left[ \mu_j g_{p,j} ({\bf u}_\infty) \right]^2 + \\
\hspace{20mm}\displaystyle \sum\limits_{i = 1}^{n_u} \left[ ( \zeta^L_i (u^L_i - u_{\infty,i}) )^2 + ( \zeta^U_i (u_{\infty,i} - u^U_i) )^2  \right] \Big ) \\
= \displaystyle \sum\limits_{j = 1}^{n_g} \left[ \tilde \mu_j g_{p,j} ({\bf u}_\infty) \right]^2 + \\
\hspace{20mm}\displaystyle \sum\limits_{i = 1}^{n_u} \left[ (\tilde \zeta^L_i (u^L_i - u_{\infty,i}) )^2 + (\tilde \zeta^U_i (u_{\infty,i} - u^U_i) )^2  \right] \\
= \displaystyle \sum\limits_{j = 1}^{n_g} \left[ \tilde \mu_j g_{p,j} ({\bf u}_\infty) \right]^2 \leq \displaystyle \sum\limits_{j = 1}^{n_g} \left( \tilde \mu_j \epsilon_j \right)^2
\end{array},
\end{equation}

where the steps may be justified as follows:

\begin{itemize}
\item Setting the Lagrange multipliers equal to their $(\tilde \cdot)$ analogues leads to taking an infimum over a set with tighter constraints, which cannot lower the value of the infimum.
\item The infimum may be evaluated by substituting in the $(\tilde \cdot)$ values and noting that the gradient of the Lagrangian is ${\bf 0}$ for this choice, leaving only the complementary slackness terms.
\item All of the terms corresponding to the box constraints may be removed since the constraints are either active with a value of 0 or have $\tilde \zeta$ coefficients equal to 0 by definition.
\item As all of the $\epsilon$-active constraints must have values that are smaller, in absolute value, than the corresponding ${\boldsymbol \epsilon}$ values, the bound may be restated further in terms of ${\boldsymbol \epsilon}$ (the $\epsilon$-inactive constraints have corresponding $\tilde \mu$ values of 0 by definition).
\end{itemize}

Clearly, this upper bound approaches 0 as ${\boldsymbol \epsilon} \rightarrow {\bf 0}$ provided that the relevant $\tilde \mu_j$ are finite for any realization of the algorithm. This is guaranteed by the finite nature of (\ref{eq:lagana}), which may only have infinite coefficients for those uncertain constraints with a zero gradient. However, this has already been ruled out by Assumption A3 as ${\boldsymbol \epsilon} \rightarrow {\bf 0}$. Since the upper bound in (\ref{eq:KKTerrup}) tends to 0, it follows that the error in (\ref{eq:KKTerr}) does as well, thereby proving (\ref{eq:KKTerr0}). \hfill $\square$      

\end{enumerate}

\end{document}